\documentclass[final,3p,times]{elsarticle}
\usepackage{amssymb}
\usepackage[reqno]{amsmath}
\usepackage{amsfonts}
\usepackage{amsthm}
\usepackage{graphicx}
\usepackage{psfrag}
\usepackage[notcite,notref]{showkeys}
\usepackage{showkeys}
\newcommand{\Sp}{Sp(2n)}

\newcommand{\bnn}{(0\,\,I)^T}

\newcommand{\tQ}{\tilde Q}

\newcommand{\bR}{{\mathbb R}}

\newcommand{\cA}{{\mathcal A}}
\newcommand{\cB}{{\mathcal B}}
\newcommand{\cC}{{\mathcal C}}
\newcommand{\cD}{{\mathcal D}}
\newcommand{\cH}{{\mathcal H}}

\newcommand{\cS}{{\mathcal S}}

\newcommand{\cJ}{{\mathcal J}}

\newcommand{\hQ}{{\hat Q}}
\newcommand{\hX}{{\hat X}}
\newcommand{\tX}{{\tilde X}}
\newcommand{\hU}{{\hat U}}
\newcommand{\hY}{{\hat Y}}

\newcommand{\hZ}{{\hat Z}}

\newcommand{\rank}{\operatorname{\text{\rm rank}}}
\newcommand{\sign}{\operatorname{\text{\rm sign}}}
\newcommand{\Ker}{\operatorname{\text{\rm Ker}}}
\newcommand{\const}{\mathop{\mathrm{const}}\nolimits}
\newcommand{\IM}{\operatorname{\text{\rm Im}}}

\newcommand{\ind}{\operatorname{\text{\rm ind}}}
\newcommand{\diag}{\operatorname{\text{\rm diag}}}

\newcommand{\smat}[4]{\left
(\begin{smallmatrix}#1&#2\\#3&#4\end{smallmatrix}\right)}

\newtheorem{theorem}{Theorem}[section]
\newtheorem{corollary}[theorem]{Corollary}
\newtheorem{lemma}[theorem]{Lemma}
\newtheorem{remark}[theorem]{Remark}
\newtheorem{definition}[theorem]{Definition}
\newtheorem{example}[theorem]{Example}
\newtheorem{proposition}[theorem]{Proposition}
\numberwithin{equation}{section}

\begin{document}

\begin{frontmatter}

\title{Cyclic sums of comparative indices and their applications}

\author{Julia Elyseeva}
\address{Department of Applied Mathematics, Moscow State University of Technology, Vadkovskii per. 3a, 101472, Moscow, Russia}

\address{{Department of Mathematics and Statistics, Faculty of Science },{Masaryk University}, Kotl\'{a}\v{r}sk\'{a} 2, CZ-61137 {{Brno},{Czech Republic}}\\e-mail:{elyseeva@gmail.com}}

\begin{abstract}
In this paper we generalize the notion of the comparative index for the pair of Lagrangian subspaces which  has fundamental applications in oscillation theory of symplectic difference systems and  linear differential Hamiltonian  systems. We introduce cyclic sums $\mu_c^{\pm}(Y_1,Y_2,\dots,Y_m),\,m\ge 2$ of the comparative indices for the set of $n-$ dimensional Lagrangian subspaces. We formulate and prove main properties of the cyclic sums, in particular,  we state connections of the cyclic sums  with the Kashiwara index. The main results of the paper connect the cyclic sums of the comparative indices with the number of positive and negative eigenvalues of $mn\times mn$   symmetric matrices defined in terms of the Wronskians $Y_i^T\cJ Y_j,$  $i,j=1,\dots,m.$ We also present first applications of the cyclic sums of the comparative indices in the oscillation theory of the  discrete symplectic systems  connecting the number of focal points of their principal solutions with the negative and positive inertia of symmetric matrices.
\end{abstract}
\end{frontmatter}
\section{Introduction}\label{sec0}
In this paper we develop an important concept from the
matrix analysis, the comparative index \cite{jE07,Elyseeva1}, \cite[chapter 3]{DEH} which has fundamental applications in the oscillation and spectral theory of the symplectic difference systems
\begin{equation}\label{SDS}
  y_{k+1}=\cS_k y_k,\quad y_k\in \mathbb{R}^{2n},\quad\cS_k\in \mathbb{R}^{2n\times 2n},\quad \cS_k^T\cJ\cS_k=\cJ,\,k=0,\dots,N,\quad \;\cJ=\begin{pmatrix}
                                  0 & I \\
                                  -I & 0 \\
                                \end{pmatrix}
\end{equation}
 as well as in the oscillation theory of the linear differential Hamiltonian systems
\begin{equation}\label{Ham1}
     y'=\cJ\cH(t)y,\,t\in [a,b],\quad y=\binom{x(t)}{u(t)},\quad \cH(t)=\cH^T(t)
\end{equation}
which are continuous counterparts of \eqref{SDS}.

In \cite{Elyseeva1}, \cite[chapter 3]{DEH} the comparative index was introduced for the pair of $2n\times n$ matrices $Y,\,\hY$ which obey the conditions
\begin{equation}\label{conj}
 w(Y,Y)=0\quad
{\rank}Y=n,\quad w(\hY,\hY)=0\quad
{\rank}\hY=n,\quad w(Y,\hY):=Y^T\cJ \hY.
\end{equation}
The matrices $ Y,\,\hY$ whose columns form bases of Lagrangian subspaces $L,\,\hat L$ in $\mathbb{R}^{2n}$ will be referred to as \textit{frames} for $L_1,\,L_2$. They can also be regarded as \textit{conjoined}  bases of \eqref{SDS} and \eqref{Ham1} and the matrix $w(Y,\hY)$ in \eqref{conj} is called the \textit{Wronskian} of $Y$ and $\hY$ according to the terminology from the oscillation theory of \eqref{SDS} and \eqref{Ham1} (see \cite{DEH} and \cite{Reid}).  The advantage of the comparative index lies in the fact that it allows to derive classical separation and comparison results for conjoined bases of \eqref{SDS} and  \eqref{Ham1} in the form of explicit relations between the multiplicities of their focal points, see \cite{Elyseeva1,Elyseeva4,jE20,Sepitka2,Sepitka3} and \cite[Chapter 4]{DEH}. Further applications of the comparative index can be found in the spectral theory of \eqref{SDS} and \eqref{Ham1}, see \cite[Chapters 5,6]{DEH}, \cite{Elyseevaaml2019} and the reference given therein. In the recent publication \cite{pS.rSH2020?e} the notion of the comparative index was connected with the traditional \textit{Lidskii angles} \cite{Lids} for symplectic matrices, in \cite{ESSmaslov} we use the comparative index defining  the so-called \textit{oscillation numbers} (see \cite{Elyseevaaml2019,jE20}) for continuous Lagrangian paths and connect the oscillation numbers with the Maslov index in \cite{pH.yL.aS2017}.

According to \cite{Elyseeva1}, we define the
{{comparative} index} for $Y,\,\hY$ partitioned into $n\times n$ blocks according to $Y=\binom{X}{U},$ $\hY=\binom{\hX}{\hU}$  using the notation
\begin{equation}
\label{compind1}
    \mathcal M=(I-X^{\dag}X)w(Y,\hY),\quad
     \mathcal T=I-\mathcal M^{\dag}\mathcal M,\quad
    \mathcal P =\mathcal T(w^T(Y,\hY)X^{\dag}\hX)\mathcal T,
\end{equation}
where the dagger $\dag$ stands for the Moore-Penrose pseudoinverse of matrices, see \cite{BenIasr}.
The {comparative index} and the \textit{dual comparative index} are defined as follows
 \begin{equation}\label{defmu}
   \mu(Y,\hat Y)=\mu_1(Y,\hat
Y)+\mu_2(Y,\hat Y),\,\mu_1(Y,\hat Y)={\rm{rank}}\mathcal M,\quad \mu_2(Y,\hat Y)={\rm{ind}}\mathcal P,
 \end{equation}
\begin{equation}\label{defmu*}
   \mu^*(Y,\hat Y)=\mu_1(Y,\hat
Y)+\mu_2^*(Y,\hat Y),\,\quad \mu_2^*(Y,\hat Y)={\rm{ind}}(-\mathcal P),
 \end{equation}
 where $\ind A$ denotes the index, i.e., the number of negative eigenvalues of the symmetric matrix $A=A^T.$

In this paper we introduce the cyclic sums of the comparative indices and the dual comparative indices for $m$ matrices $Y_1,Y_2,\dots,Y_m,\,m\ge 2$ with condition \eqref{conj} according to the following definition.
\begin{definition}\label{cyclsumd}
  Consider  $2n\times n$ matrices $Y_k,\,k=1,2,\dots,m,\,m\ge 2$ which obey condition \eqref{conj}. We define the cyclic sums (of the first kind) as follows
\begin{equation}\label{cyclsum}\begin{aligned}
 \mu_c^-(Y_1,Y_2,\dots,Y_m):= \mu(Y_1,Y_2)+\mu(Y_2,Y_3)+\dots+\mu(Y_{m-1},Y_m)+\mu(Y_m,Y_1),\\
 \mu_c^+(Y_1,Y_2,\dots,Y_m):= \mu^*(Y_1,Y_2)+\mu^*(Y_2,Y_3)+\dots+\mu^*(Y_{m-1},Y_m)+\mu^*(Y_m,Y_1).
\end{aligned}\end{equation}
By a similar way we introduce the cyclic sums (of the second kind)
\begin{equation}\label{cyclsummod}\begin{aligned}
 \nu_c^-(Y_1,Y_2,\dots,Y_m):= \mu(Y_1,Y_2)+\mu(Y_2,Y_3)+\dots+\mu(Y_{m-1},Y_m)-\mu(Y_1,Y_m),\\
 \nu_c^+(Y_1,Y_2,\dots,Y_m):= \mu^*(Y_1,Y_2)+\mu^*(Y_2,Y_3)+\dots+\mu^*(Y_{m-1},Y_m)-\mu^*(Y_1,Y_m).
\end{aligned}\end{equation}
\end{definition}
The first  property of \eqref{cyclsum} for $m=2$  can be found among the  properties of the comparative index (see \cite[p.448, property 5]{Elyseeva1},\cite[Theorem~3.5(v)]{DEH})
\begin{equation}\label{n=2}\begin{aligned}
 \mu_c^{-}(Y_1,Y_2)=\mu(Y_1,Y_2)+\mu(Y_2,Y_1)&=\mu^{+}_c(Y_2,Y_1)= \mu^*(Y_2,Y_1)+\mu^*(Y_1,Y_2)\\&=\rank w(Y_1,Y_2)=\ind \smat{0}{w(Y_1,Y_2)}{w^T(Y_1,Y_2)}{0},
\end{aligned}\end{equation}
where the last equality follows from \cite[Fact 5.8.20]{Bernstein}. Observe that by \eqref{n=2} we have the main connection for the cyclic sums \eqref{cyclsum} and \eqref{cyclsummod}
\begin{equation}\label{add5}
  \mu_c^{\pm}(Y_1,Y_2,\dots,Y_m)=\nu_c^{\pm}(Y_1,Y_2,\dots,Y_m)+\rank w(Y_m,Y_1),
\end{equation}
in particular, $\nu_c^{\pm}(Y_1,Y_2)=0.$

For $m=3$ we prove that the cyclic sums $\mu_c^{+}(Y_1,Y_2,Y_3)$ and $\mu_c^{-}(Y_1,Y_2,Y_3)$ defined by \eqref{cyclsum} coincide with the \textit{positive} $i_{+}(B(x,x))$ and \textit{negative} $i_{-}(B(x,x))$ inertia of  the the quadratic form
\begin{equation}\label{bilin}
 B(x,x):=B((x_1,x_2,x_3), (x_1,x_2,x_3))=w(x_1,x_2)+w(x_2,x_3)+w(x_3,x_1),\quad w(x_i,x_j)=x_i^T\cJ x_j
\end{equation}
defined on $(x_1,x_2,x_3)\in L_1\oplus L_2\oplus L_3,$ where $L_1,\,L_2,\,L_3$ are Lagrangian subspaces in $\mathbb{R}^{2n}$ with the frames $Y_1,\,Y_2,\,Y_3.$
Recall that \textit{the Kashiwara index} \cite[Definition A.3.1]{Kash} is defined as the signature
\begin{equation}\label{kash}
  \tau(L_1,L_2,L_3)=\sign(B(x,x))=i_{+}(B(x,x))-i_{-}(B(x,x))
  \end{equation}
of the quadratic form $B(x,x)$, i.e.,  the difference of positive and negative eigenvalues of $B(x,x)$ and therefore $\tau(L_1,L_2,L_3)$ can be presented in terms of the cyclic sums $\mu_c^{\pm}(Y_1,Y_2,Y_3)$
\begin{equation}\label{kashcycls3}
 \tau(L_1,L_2,L_3)= \mu_c^{+}(Y_1,Y_2,Y_3)-\mu_c^{-}(Y_1,Y_2,Y_3).
\end{equation}
We prove that \eqref{kashcycls3} is valid for the general case  $m\ge 3$
\begin{equation}\label{kashcyclsm}
 \tau(L_1,L_2,\dots,L_m)= \mu_c^{+}(Y_1,Y_2,\dots,Y_m)-\mu_c^{-}(Y_1,Y_2,\dots,Y_m),
\end{equation}
where  $\tau(L_1,L_2,\dots,L_m)$ is defined by the formula (see \cite[Definition A.3.7]{Kash})
\begin{equation}\label{kashm}
  \tau(L_1,L_2,\dots,L_m)=\sum\limits_{j=2}^{m-1}\tau(L_1,L_j,L_{j+1}).
\end{equation}

The main result of the paper connects $\mu_c^{\pm}(Y_1,Y_2,\dots,Y_m)$ with the  negative and positive inertia $i_{\pm}$ of a $(m n)\times (mn)$ symmetric matrix according to the following theorem.

\begin{theorem}\label{maincon}
Define the $(m n)\times (mn)$ symmetric matrix
\begin{equation}\label{Sm}
  S_{1,2,\dots,m}=\begin{pmatrix}
                    0 & w_{1,2} & w_{1,3} & \dots & w_{1,m} \\
                    w^T_{1,2} & 0 & w_{2,3} & \dots & w_{2,m} \\
                    \dots & \dots & \dots & \dots & \dots \\
                    w^T_{m-1,1} & w^T_{m-1,2} & \dots & 0 & w_{m-1,m} \\
                    w^T_{m,1} & w^T_{m,2} & \dots & w^T_{m,m-1} & 0 \\
                  \end{pmatrix},\quad w_{i,j}:=w(Y_i,Y_j)=Y_i^T\cJ Y_j,
\end{equation}
then for cyclic sums \eqref{cyclsum} we have
\begin{equation}\label{Scyclm}
  \mu_c^{\pm}(Y_1,Y_2,\dots,Y_m)=i_{\mp}(S_{1,2,\dots,m})=\ind(\pm \;S_{1,2,\dots,m}).
\end{equation}
\end{theorem}

Formula \eqref{Scyclm} presents fundamental connections between the comparative index theory and matrix linear algebra. The proof is based on the main theorem of the comparative index (see \cite[Theorem 2.2]{Elyseeva1}) which implies that cyclic sums \eqref{cyclsum}, \eqref{cyclsummod}  are  invariant  with respect to arbitrary
symplectic  transformations, compare with \cite[Property III (Symplectic invariance)]{seC.rL.eyM1994}. A similar index result is also proved for the cyclic sums of the second kind (see Theorem~\ref{maincon23}), in the proof we apply to \eqref{Sm} the results by Y. Tian (see \cite[Theorem 2.3]{Tian})
concerning evaluations of inertias of block symmetric matrices. Theorems~\ref{maincon},~\ref{maincon23}  can be used as a new tool for  computations  of  cyclic sums \eqref{cyclsum}, \eqref{cyclsummod} and  the Kashiwara  indices defined by \eqref{kash} and \eqref{kashm}.

The Kashiwara  and the \textit{H\"{o}rmander} indices \cite{Horm}  are traditionally related to the Maslov index theory for Lagrangian paths, see \cite{Kash,seC.rL.eyM1994,bBB.cZ2018,pH.yL.aS2017,Howardnew,Zhou}.  We conjecture that the cyclic sums \eqref{cyclsum}, \eqref{cyclsummod} and their properties proved in this paper will be a useful complement to the theory of the Maslov index and it's applications which are the subject of our future investigations.

In the present paper we  concentrate on  new applications of the cyclic sums \eqref{cyclsum}, \eqref{cyclsummod} in the oscillation theory of \eqref{SDS}. For \eqref{SDS} we consider arbitrary fundamental symplectic matrices $\mathcal Z_k\in\Sp,\,k=0,1,\dots,N+1$ and show that the cyclic sums
\begin{equation}\label{secka}
  \mu_c^{-}(\mathcal Z_0^{-1}\bnn,\mathcal Z_1^{-1}\bnn,\dots,\mathcal Z_{N+1}^{-1}\bnn)=\mu_c^{+}(\mathcal Z_{N+1}^{-1}\bnn,\mathcal Z_N^{-1}\bnn,\dots,\mathcal Z_{0}^{-1}\bnn)
\end{equation}
and
\begin{equation}\label{seckb}
  \nu_c^{-}(\mathcal Z_0^{-1}\bnn,\mathcal Z_1^{-1}\bnn,\dots,\mathcal Z_{N+1}^{-1}\bnn)=\nu_c^{+}(\mathcal Z_{N+1}^{-1}\bnn,\mathcal Z_N^{-1}\bnn,\dots,\mathcal Z_{0}^{-1}\bnn)
\end{equation}
are invariant with respect to $\mathcal Z_k\in\Sp,\,k=0,1,\dots,N+1$ and present the maximal and minimal numbers of focal points of conjoined bases of \eqref{SDS}. These numbers also coincide with the numbers of (forward) focal points of the principal solutions of \eqref{SDS} at $k=N+1$ and  $k=0,$ respectively (see Theorem~\ref{cycldisc}). Observe that the minimal number of focal points described by the cyclic sums of the second kind \eqref{seckb} is highly important in the oscillation and spectral theory of \eqref{SDS}, see e.g. \cite{BDosly,BDosly1,Kratz1,mB.wK.rSH12} and \cite[Chapters 4,5]{DEH}. In particular, by the \textit{Reid Roundabout Theorem}, system \eqref{SDS}  is \textit{disconjugate} on $[0,N+1]$ if and only if the number of (forward) focal points of the principal solution at $k=0$ is equal to zero (see \cite[Theorem 1]{BDosly},\cite[Theorem 2.36]{DEH}). In \cite{BDosly1,Hilsch5} this result was interpreted as nonnegative definiteness of  $n(N+1)\times n(N+1)$ symmetric matrices associated with \eqref{SDS}. In this paper, applying the index results for cyclic sums (Theorems~\ref{maincon},~\ref{maincon23}) we present a  generalization of the results in \cite{BDosly1,Hilsch5} connecting the number of (forward) focal points of the principal solution at $k=0$ with the index of  $n(N+1)\times n(N+1)$ symmetric matrices, see Theorem~\ref{indfocd}.

The paper is organized as follows. In the next section we recall in more details the notion of the comparative index and its properties,  formulate and prove main properties of cyclic sums \eqref{cyclsum}, \eqref{cyclsummod} (see Propositions~\ref{invar},~\ref{invar1},~\ref{recrel}). In Section~\ref{sec2} we present the proof of Theorem~\ref{maincon}, formulate and prove Theorem~\ref{maincon23} connecting the cyclic sums \eqref{cyclsummod} with the index of symmetric matrices. In Section~\ref{sec3} we present applications of the cyclic sums, in particular, we prove connections \eqref{kashcycls3}, \eqref{kashcyclsm} and similar connections for the cyclic sums of the second kind (see Theorem~\ref{kashcycl}). In  Section~\ref{sec3} we also present the above applications to the oscillation theory of \eqref{SDS}, in particular, we prove Theorems~\ref{cycldisc},~\ref{indfocd} and their corollaries.
\section{Properties of cyclic sums}\label{sec1}
We will use the following notation. For a matrix $A,$ we denote by
$A^T,\,A^{-1},$  $A^{\dag}, \rank A,$ $\,\Ker A,\,\IM A,$
$ \ind A,$ $\sign A,$ $ A \ge~0, A \le 0,$ respectively, its transpose, inverse,
 Moore-Penrose pseudoinverse,  rank (i.e.,
the dimension of its image), kernel, image, index (i.e., the number
of its negative eigenvalues), signature (i.e., the difference between the numbers of positive and negative eigenvalues of $A$),  positive semidefiniteness, negative
semidefiniteness and we use the notation $E_A=I-AA^{\dag},\,F_A=I-A^{\dag}A$ for the orthogonal projectors on to the null spaces of $A^T$ and $A$, respectively.
We also use the notation $\Sp$ for the real matrix symplectic group in dimension $2n$ and $\Delta A_k=A_{k+1}-A_k$ for the forward differences.

Firstly we recall the definition of the  comparative index in more details, see \cite{Elyseeva1}. There are four equivalent definitions of the comparative index according to \cite[Theorem 2.1]{Elyseeva1}. Here we mention the definition in terms of symmetric solutions $Q,\;\hQ$  of the matrix equations
 $X^TQX=X^TU,\quad \hX^T\hQ\hX=\hX^T\hU.$
Introduce the matrices
\begin{equation}
\label{compind}
    \tilde {\mathcal M}=(I-XX^{\dag})\hat X,\quad
     \mathcal T=I-\tilde{\mathcal M}^{\dag}\tilde{\mathcal M},\quad
    {\mathcal P} =\mathcal T\hX^T(\hQ-Q)\hX \mathcal T.
\end{equation}
Then $\mu(Y,\hY)$ and $\mu^*(Y,\hY)$ can be defined according to \eqref{defmu}, \eqref{defmu*}, where $\mathcal M$ in \eqref{compind1} is replaced by $\tilde {\mathcal M}$  (observe that the matrices $\mathcal T$, $\mathcal P$ stay the same as in \eqref{compind1}  after the replacement of $\mathcal M$ by $\tilde {\mathcal M}$, see \cite[Theorem 2.1]{Elyseeva1}). Moreover, the matrix $\mathcal P$ in \eqref{compind1} has the equivalent representation given by \eqref{compind}.

Applying the formula (see \cite{Styan})
\begin{equation}\label{stian1}\begin{aligned}
  \rank (A\;B)&=\rank A + \rank (E_A B),\,A\in \mathbb{C}^{l\times n},\,B\in \mathbb{C}^{l\times k}
 \end{aligned}\end{equation}
one can define the first addends $\mu_1(Y,\hY)$ in \eqref{defmu}, \eqref{defmu*} as follows
\begin{equation}\label{newmu1}\begin{aligned}
  \mu_1(Y,\hY)&=\rank  \tilde {\mathcal M}=\rank (X\;\hX)-\rank (X)=\rank   {\mathcal M}=\rank (X^T\;w(Y,\hY))-\rank (X).
\end{aligned}\end{equation}

Incorporating formulas for the inertia of block symmetric matrices with the blocks $A=A^*\in \mathbb{C}^{k\times k},$ $D=D^*\in \mathbb{C}^{l\times l},$ and $B\in \mathbb{C}^{k\times l}$
\begin{equation}\label{tianin}\begin{aligned}
  i_{\pm} \begin{pmatrix}
         A & B \\
         B^T & D \\
       \end{pmatrix}&=i_{\pm} (A)+i_{\pm}\begin{pmatrix}
                                    0 & M \\
                                    M^T & D-B^T A^{\dag} B \\
                                  \end{pmatrix},\\i_{\pm}\begin{pmatrix}
                                    0 & M \\
                                    M^T & D-B^T A^{\dag} B \\
                                  \end{pmatrix}&=\rank (M)+i_{\pm} (F_M (D-B^T A^{\dag} B)F_M),\,\,M=E_A B,
\end{aligned}\end{equation}
(see \cite[Theorem 2.3]{Tian} and \cite{Elyseeva1}, where \eqref{tianin}   was derived for the proof of the comparative index properties for the special case $k=l=n$) the comparative index $\mu(Y,\hY)$ and the dual comparative index $\mu^*(Y,\hY)$ can be shortly defined as follows (see \cite[Lemma 3.14]{DEH})
\begin{equation}\label{shortdefc}
  \mu(Y,\hY)=i_{-}\begin{pmatrix}
                    0 & \tilde {\mathcal M}\\
                    {\tilde {\mathcal M}}^T & \hX^T(\hQ-Q)\hX \\
                  \end{pmatrix},\quad \mu^*(Y,\hY)=i_{+}\begin{pmatrix}
                    0 & \tilde {\mathcal M}\\
                    {\tilde {\mathcal M}}^T & \hX^T(\hQ-Q)\hX \\
                  \end{pmatrix},
\end{equation}
where $\tilde {\mathcal M}$ is given by \eqref{compind} ( $\tilde {\mathcal M}$ can be replaced by $\mathcal M$ in \eqref{compind1}). In particular, for the case $\det X~\ne~0,\;\det \hX~\ne~0$ we have
\begin{equation}\label{compnonsing}
 \mu(Y,\hY)=\mu_2(Y,\hY)=\mu_2^*(\hY,Y)=\mu^*(\hY,Y)=\ind (\hQ-Q),\,Q=UX^{-1},\;\hQ=\hU\hX^{-1}.
\end{equation}
Other special cases for the comparative index are the following
\begin{equation}\label{speccase}
  \mu(Y,\bnn)=0,\;\mu(\bnn,\hY)=\mu_1(\bnn,\hY)=\rank \hX,\;\mu(\cJ\bnn,\hY)=\mu_2(\cJ\bnn,\hY)=\ind(\hX^T\hU).
\end{equation}
For the convenience we collect some properties of the comparative index which we will use in the subsequent proofs (see  \cite[p.448]{Elyseeva1} and \cite[Theorem 3.5 and formula (3.34), p.165]{DEH} ).
\begin{lemma}\label{prop ind}
 For $Z,\,\hZ$  and $Y=\binom{X}{U}=Z\bnn,\,\hY=\binom{\hX}{\hU}=\hZ\bnn$  we have the following properties.
 \begin{itemize}
 \item [(i)] $  \mu_l(LYC,L\hY\hat C)=\mu_l(Y,\hY),\,l=1,2,\;\det C\ne 0,\,\det \hat C\ne 0,$ where  $L$ is an arbitrary symplectic lower block-triangular matrix.
     \item[(ii)] $\mu(Y,\hY)+\mu(\hY,Y)=\rank w(Y,\hY),\quad w(Y,\hY)=Y^T\cJ\hY,$
       \item [(iii)] $\mu(Y,\hY)=\rank \hX-\rank X+\mu^*(\hY,Y),$
   \item[(iv)]    $ \mu_l(Y,\hat Y)=\mu^{*}_l(Z^{-1}\bnn,Z^{-1}\hat
    Y),\,l=1,2,$
       \item[(v)] $\mu(Y,\hY)+\mu^*(Y,\hY)=\rank w(Y,\hY)+\rank \hX-\rank X.$
       \item[(vi)] $0 \le \mu(Y,\hY)\le \min(\rank w(Y,\hY),\rank \hX)\le n.$
  \end{itemize}
\end{lemma}

Observe  that according to a \textit{duality principle}  Lemma~\ref{prop ind} holds also for the dual index $\mu^*(Y,\hY)$ (where we use that according to the definition $(\mu^*(Y,\hY))^*=\mu(Y,\hY)$), see \cite{Elyseeva1} and \cite[Chapter 3]{DEH}.

Introduce $Z_i\in\Sp$ associated with $Y_i$ in Definition~\ref{cyclsum} according to
\begin{equation}\label{matr}
  Y_i=Z_i\bnn,\,i=1,2,\dots,m.
\end{equation}
We have the connection
\begin{equation}\label{iverseandwrons}
  (I\;0)Z_i^{-1}Y_j=-w_{i,j},\,w_{i,j}:=w(Y_i,Y_j),
\end{equation}
where we use that $Z_i^{-1}=-\cJ Z_i^T\cJ$ for $Z_i\in\Sp.$ Remark that $Z_i\in\Sp$ are not uniquely defined by \eqref{matr}. We have $Y_i=Z_i L\bnn $ for arbitrary symplectic unit lower block-triangular matrix $L=\smat{I}{0}{Q}{I},\,Q=Q^T.$ In subsequent computations we will use $Z^{-1}_i$ taking in mind that  results of the computations do not depend on the choice of $Z_i$ in \eqref{matr} by Lemma~\ref{prop ind}(i).
\begin{example}\label{ex1}
Let $Y_1,\,Y_2,\,Y_3$ be $2n\times n$ matrices with conditions \eqref{conj}.  Consider the definition of the comparative index $\mu(Z_3^{-1}Y_1,Z_3^{-1}Y_2)$ and the dual comparative index $\mu^*(Z_3^{-1}Y_1,Z_3^{-1}Y_2)$, where $Z_3$ obeys \eqref{matr}.

Putting $Y:=Z_3^{-1}Y_1,\,\hY:= Z_3^{-1}Y_2$ we have by \eqref{compind}, \eqref{compind}, and \eqref{newmu1}
\begin{equation}\label{compwron1a}\begin{aligned}
  \mu_1(Z_3^{-1}Y_1,Z_3^{-1}Y_2)&=\rank (E_{w_{3,1}}w_{3,2})=\rank (w_{3,1}\,w_{3,2})-\rank w_{3,1}\\&=\rank (F_{w_{3,1}}w_{1,2})=\rank (w_{1,3}\,w_{1,2})-\rank w_{1,3},
\end{aligned}\end{equation}
where we computed the Wronskian of $Z_3^{-1}Y_1,Z_3^{-1}Y_2$ according to  $w(Z_3^{-1}Y_1,Z_3^{-1}Y_2)=w(Y_1,Y_2).$

Next we consider the second addend given by \eqref{compind1}, \eqref{defmu}, \eqref{defmu*}
\begin{equation}\label{compwron2}\begin{aligned}
  \mu_2(Z_3^{-1}Y_1,Z_3^{-1}Y_2)&=\ind ( F_{\mathcal M}\; \mathcal D\; F_{\mathcal M}),\quad \mathcal M= F_{w_{3,1}}w_{1,2},\\\mu_2^*(Z_3^{-1}Y_1,Z_3^{-1}Y_2)&=\ind (- F_{\mathcal M}\; \mathcal D\; F_{\mathcal M}), \quad \mathcal D=w^T_{1,2}w^{\dag}_{3,1}w_{3,2}
\end{aligned}\end{equation}
and define $\mu(Z_3^{-1}Y_1,Z_3^{-1}Y_2)$ and  $\mu^*(Z_3^{-1}Y_1,Z_3^{-1}Y_2)$ given by \eqref{compwron1a} and \eqref{compwron2} according to \eqref{tianin}
\begin{equation}\label{ind1}\begin{aligned}
  \mu(Z_3^{-1}Y_1,Z_3^{-1}Y_2)&=i_{-}\begin{pmatrix}
                                      0 & \mathcal M \\
                                      \mathcal M^T & \mathcal D+ \mathcal D^T \\
                                    \end{pmatrix},\quad \mu^*(Z_3^{-1}Y_1,Z_3^{-1}Y_2)&=i_{+}\begin{pmatrix}
                                      0 & \mathcal M \\
                                      \mathcal M^T & \mathcal D+ \mathcal D^T\\
                                    \end{pmatrix},
\end{aligned}\end{equation}
where we use that $F_{\mathcal M}DF_{\mathcal M}$ is symmetric and then after the application of \eqref{tianin} we will have $F_{\mathcal M}DF_{\mathcal M}+F_{\mathcal M}D^TF_{\mathcal M}=2F_{\mathcal M}DF_{\mathcal M}.$ Remark that \eqref{ind1} can be derived directly from Theorem~\ref{maincon} for $m=3$ together with other representations of the given type incorporating properties of the cyclic sums, see Section~\ref{sec2}.
\end{example}

Next we formulate  properties of \eqref{cyclsum}, \eqref{cyclsummod} based on Lemma~\ref{prop ind}.
\begin{proposition}\label{invar1}
The cyclic sums \eqref{cyclsum}, \eqref{cyclsummod} obey the following properties.
\par{(i)} For arbitrary nonsingular matrices $C_k\in \mathbb{R}^{n\times n},\,k=1,2,\dots,m$ we have
\begin{equation}\label{inv5}\begin{aligned}
 \mu_c^{\pm}(Y_1,Y_2,\dots,Y_m)=& \mu_c^{\pm}(Y_1C_1,Y_2C_2,\dots,Y_mC_m),\quad
  \nu_c^{\pm}(Y_1,Y_2,\dots,Y_m)= \nu_c^{\pm}(Y_1C_1,Y_2C_2,\dots,Y_mC_m).
\end{aligned}\end{equation}
\par{(ii)} According to the definition in \eqref{cyclsum} we have the following invariant property with respect to  cyclic permutations
\begin{equation}\label{inv2}\begin{aligned}
   \mu_c^{\pm}(Y_1,Y_2,\dots,Y_m)= \mu_c^{\pm}(Y_m,Y_1,Y_2\dots,Y_{m-1})&=\dots =\mu_c^{\pm}(Y_2,Y_3,\dots,Y_m,Y_1).
\end{aligned}\end{equation}
\par{(iii)} We  have
\begin{equation}\label{inv3}
  \mu_c^{\pm}(Y_1,Y_2,\dots,Y_m)=\mu_c^{\mp}(Y_m,Y_{m-1},\dots,Y_1),\quad \nu_c^{\pm}(Y_1,Y_2,\dots,Y_m)=\nu_c^{\mp}(Y_m,Y_{m-1},\dots,Y_1).
\end{equation}
\par{(iv)} We have
\begin{equation}\label{con1+}\begin{aligned}
  \mu_c^{-}(Y_1,Y_2,\dots,Y_m)+\mu_c^{+}(Y_1,Y_2,\dots,Y_m)=\sum\limits_{j=1}^{m-1}\rank w(Y_j,Y_{j+1})+\rank w(Y_m,Y_1),\\\nu_c^{-}(Y_1,Y_2,\dots,Y_m)+\nu_c^{+}(Y_1,Y_2,\dots,Y_m)=\sum\limits_{j=1}^{m-1}\rank w(Y_j,Y_{j+1})-\rank w(Y_m,Y_1)
\end{aligned}\end{equation}
\end{proposition}
\begin{proof}
The proof of \eqref{inv5} follows from Lemma~\ref{prop ind}(i) and the definitions of the cyclic sums \eqref{cyclsum}, \eqref{cyclsummod}.

The proof of \eqref{inv2} follows from  \eqref{cyclsum}. Indeed, the cyclic permutations of $Y_k$ do not change the order of the pairs $Y_k,Y_{k+1}$ and $Y_m,Y_1$ in \eqref{inv2}.

For the proof of \eqref{inv3} we apply Lemma~\ref{prop ind}(iii)
\begin{equation*}
  \mu(Y_k,Y_{k+1})=\mu^*(Y_{k+1},Y_{k})+\rank X_{k+1}-\rank X_k,\,k=1,\dots,m-1,\;\mu(Y_m,Y_{1})=\mu^*(Y_1,Y_{m})+\rank X_{1}-\rank X_{m}.
\end{equation*}

Summing the previous identities  we derive
\begin{equation}\label{add15}\begin{aligned}
  \mu_c^{-}(Y_1,Y_2,\dots,Y_m)&=\sum\limits_{k=1}^{m-1} \mu(Y_k,Y_{k+1}) + \mu(Y_m,Y_1)=\sum\limits_{k=1}^{m-1} \mu^*(Y_{k+1},Y_{k}) + \mu^*(Y_1,Y_m)+s_m=\mu^{+}(Y_m,Y_{m-1},\dots,Y_1)+s_m,
\end{aligned}\end{equation}
where
  $s_m=\sum\limits_{k=1}^{m-1}\Delta(\rank X_k)+\rank X_{1}-\rank X_{m}=0.$
The proof of   equality \eqref{inv3} for $\mu_c^{+}(Y_1,Y_2,\dots,Y_m)$ based on a dual form of Lemma~\ref{prop ind}(iii)  is similar. For $\nu_c^{\pm}(Y_1,Y_2,\dots,Y_m)$ we have according to \eqref{add5}
\begin{equation*}
  \nu_c^{\pm}(Y_1,Y_2,\dots,Y_m)=\mu_c^{\pm}(Y_1,Y_2,\dots,Y_m)-\rank w(Y_m,Y_1)=\mu_c^{\mp}(Y_m,Y_{m-1},\dots,Y_1)-\rank w(Y_1,Y_m)=\nu_c^{\mp}(Y_1,Y_2,\dots,Y_m),
\end{equation*}
where we used the first equality in \eqref{inv3} proved above.

The proof of \eqref{con1+}  is based on Lemma~\ref{prop ind}(v). We have
\begin{equation*}\begin{aligned}
  \mu(Y_k,Y_{k+1})+\mu^*(Y_{k},Y_{k+1})&=\rank X_{k+1}-\rank X_k+\rank w(Y_k,Y_{k+1}),\,k=1,\dots,m-1,\\\mu(Y_m,Y_{1})+\mu^*(Y_m,Y_{1})&=\rank X_1-\rank X_m+\rank w(Y_m,Y_1).
\end{aligned}\end{equation*}

Summing the previous identities  we derive
\begin{equation*}\begin{aligned}
  \mu_c^{-}(Y_1,Y_2,\dots,Y_m)+\mu_c^{+}(Y_1,Y_2,\dots,Y_m)=\sum\limits_{j=1}^{m-1}\rank w(Y_j,Y_{j+1})+\rank w(Y_m,Y_1)+s_m,
\end{aligned}\end{equation*}
where $s_m=0$ is the same as in \eqref{add15}. So we proved the first identity in \eqref{con1+}. Applying \eqref{add5} we also prove the second one.
The proof is completed.
\end{proof}

As it was mentioned in Section~\ref{sec0} the main property of the cyclic sums \eqref{cyclsum} and \eqref{cyclsummod} is their invariance with respect to arbitrary symplectic transformations. The proof of the invariance is based on the main theorem of the comparative index theory (see \cite[Theorem 2.2, formulas (2.14), (2.15)]{Elyseeva1}, \cite[Theorem 3.5, Corollary 3.12, formulas (3.17),(3.26)]{DEH}).
\begin{theorem}\label{maincompind}
For arbitrary $W\in\Sp$ and $2n\times n$ matrices $Y,\;\hat Y$ with
condition \eqref{conj} we have
\begin{equation}\label{mainth}\begin{aligned}
  \mu(WY,W\hY)&=\mu(Y,\hY)+\mu(\hY,W^{-1}\bnn)-\mu(Y,W^{-1}\bnn),\\
  \mu^*(WY,W\hY)&=\mu^*(Y,\hY)+\mu^*(\hY,W^{-1}\bnn)-\mu^*(Y,W^{-1}\bnn).
\end{aligned}\end{equation}
\end{theorem}

Based on Theorem~\ref{maincompind}  we prove the following result.
\begin{proposition}[Symplectic invariance]\label{invar}
The cyclic sums \eqref{cyclsum} and \eqref{cyclsummod}
are invariant with respect to an arbitrary symplectic transformation, i.e., for arbitrary  matrix $R\in\Sp$
\begin{equation}\label{inv1}\begin{aligned}
  \mu_c^{\pm}(Y_1,Y_2,\dots,Y_m)=\mu_c^{\pm}(R^{-1}Y_1,R^{-1}Y_2,\dots,R^{-1}Y_m),\quad
  \nu_c^{\pm}(Y_1,Y_2,\dots,Y_m)=\nu_c^{\pm}(R^{-1}Y_1,R^{-1}Y_2,\dots,R^{-1}Y_m).
\end{aligned}\end{equation}
\end{proposition}
\begin{proof}
  According to Theorem~\ref{maincompind}  we have for $W:=R^{-1}$
\begin{equation}\label{ad2}\begin{aligned}
  \mu(R^{-1}Y_k,R^{-1}Y_{k+1})&=\mu(Y_k,Y_{k+1})+\mu(Y_{k+1},R\bnn)-\mu(Y_{k},R\bnn),\,k=1,\dots,m-1,\\
  \mu(R^{-1}Y_m,R^{-1}Y_{1})&=\mu(Y_m,Y_{1})+\mu(Y_{1},R\bnn)-\mu(Y_{m},R\bnn).
\end{aligned}\end{equation}
Summing the first identities for all $k=1,\dots,m-1$ and the second one we derive
\begin{equation*}
  \mu_c^{-}(R^{-1}Y_1,R^{-1}Y_2,\dots,R^{-1}Y_m)=\sum\limits_{k=1}^{m-1} \mu(R^{-1}Y_k,R^{-1}Y_{k+1}) + \mu(R^{-1}Y_m,R^{-1}Y_1)=\mu_c^{-}(Y_1,Y_2,\dots,Y_m)+s_m,
\end{equation*}
where
\begin{equation*}\begin{aligned}
  s_m=\sum\limits_{j=1}^{m-1}\Delta(\mu(Y_{j},R\bnn))+\mu(Y_{1},R\bnn)-\mu(Y_{m},R\bnn)=0.
\end{aligned}\end{equation*}
The proof for  $\mu_c^{+}(Y_1,Y_2,\dots,Y_m)$ based on the second identity in \eqref{mainth} for the dual indices is similar. The invariance of the sums $\nu_c^{+}(Y_1,Y_2,\dots,Y_m)$ follows from \eqref{add5}. Indeed, we have proved the invariance of $\mu_c^{\pm}(Y_1,Y_2,\dots,Y_m),$ then the sums $\nu_c^{\pm}(Y_1,Y_2,\dots,Y_m)=\mu_c^{\pm}(Y_1,Y_2,\dots,Y_m)-w(Y_1,Y_m)$ are invariant because of the obvious invariant property of the Wronskian $$w(R^{-1}Y_1,R^{-1}Y_m)=Y_1^TR^{-1\,T}\cJ R^{-1}Y_m=Y_1^T\cJ Y_m=w(Y_1,Y_m).$$  The proof is completed.
\end{proof}
\begin{corollary}\label{impinv}
Putting $R:=Z_m$ in Proposition~\ref{invar} we derive
the following representations for cyclic sums \eqref{cyclsum}, \eqref{cyclsummod}
\begin{equation}\label{partcase2}\begin{aligned}
  \mu_c^{-}(Y_1,Y_2,\dots,Y_m)&=\mu_c^{-}(Z_m^{-1}Y_1,Z_m^{-1}Y_2,\dots,Z_m^{-1}Y_{m-1},\bnn)=\sum\limits_{j=1} ^{m-2}\mu(Z_m^{-1}Y_j,Z_m^{-1}Y_{j+1})+\rank w(Y_m,Y_{1}),\\
  \nu_c^{-}(Y_1,Y_2,\dots,Y_m)&=\mu_c^{-}(Z_m^{-1}Y_1,Z_m^{-1}Y_2,\dots,Z_m^{-1}Y_{m-1},\bnn)=\sum\limits_{j=1} ^{m-2}\mu(Z_m^{-1}Y_j,Z_m^{-1}Y_{j+1}).
\end{aligned}\end{equation}
Similar equalities hold for $\mu_c^{+}(Y_1,Y_2,\dots,Y_m)$ and $\nu_c^{+}(Y_1,Y_2,\dots,Y_m)$ with the  dual comparative indices $\mu^*(Z_m^{-1}Y_j,Z_m^{-1}Y_{j+1})$ in the right-hand sides of   \eqref{partcase2},
in particular, it follows from the representations for $\nu_c^{\pm}(Y_1,Y_2,\dots,Y_m)$
\begin{equation}\label{nonneg}
  \nu_c^{\pm}(Y_1,Y_2,\dots,Y_m)\ge 0.
\end{equation}
\end{corollary}
\begin{remark}
\par{(i)} We have proved Proposition~\ref{invar} applying Theorem~\ref{maincompind}. Indeed the invariant property \eqref{inv1} for $\nu^{\pm}_c(Y_1,Y_2,Y_3)$ and Theorem~\ref{maincompind} are equivalent. We have for $m=3$
\begin{equation}\label{invmainTh}
  \nu^{\pm}_c(Y,\hY,W^{-1}\bnn)=\mu(Y,\hY)+\mu(\hY,W^{-1}\bnn)-\mu(Y,W^{-1}\bnn)=\nu^{\pm}_c(WY,W\hY,\bnn)=\mu(WY,W\hY),
\end{equation}
where we applied \eqref{inv1}   for $R:=W^{-1}.$
\par{(ii)} For the case $m=3$ we derive from Corollary~\ref{impinv}
\begin{equation}\label{m=3invara}\begin{aligned}
  \mu^{-}_c(Y_1,Y_2,Y_3)&=\mu^{-}_c(Z_3^{-1}Y_1,Z_3^{-1}Y_2,\bnn)=\rank w(Y_1,Y_3)+\mu(Z_3^{-1}Y_1,Z_3^{-1}Y_2),\\ \mu^{+}_c(Y_1,Y_2,Y_3)&=\mu^{+}_c(Z_3^{-1}Y_1,Z_3^{-1}Y_2,\bnn)=\rank w(Y_1,Y_3)+\mu^*(Z_3^{-1}Y_1,Z_3^{-1}Y_2),
 \end{aligned}\end{equation}
\begin{equation}\label{m=3invarb}\begin{aligned}
  \nu^{-}_c(Y_1,Y_2,Y_3)&=\nu^{-}_c(Z_3^{-1}Y_1,Z_3^{-1}Y_2,\bnn)=\mu(Z_3^{-1}Y_1,Z_3^{-1}Y_2),\\
  \nu^{+}_c(Y_1,Y_2,Y_3)&=\nu^{+}_c(Z_3^{-1}Y_1,Z_3^{-1}Y_2,\bnn)=\mu^*(Z_3^{-1}Y_1,Z_3^{-1}Y_2),
 \end{aligned}\end{equation}
 Recall that we already computed $\mu(Z_3^{-1}Y_1,Z_3^{-1}Y_2),$ $\mu^*(Z_3^{-1}Y_1,Z_3^{-1}Y_2)$ in Example~\ref{ex1}.
\end{remark}

Next we present  recurrent relations for the cyclic sums \eqref{cyclsum}.
and \eqref{cyclsummod}.

\begin{proposition}\label{recrel}
\par{(i)} For cyclic sums  \eqref{cyclsummod} we have for  any $2\le l<m$
\begin{equation}\label{rec2mod}\begin{aligned}
  \nu_c^{\pm}(Y_1,Y_2,\dots,Y_m)&=\nu_c^{\pm}(Y_1,Y_2,\dots,Y_{l})+\nu_c^{\pm}(Y_1,Y_l,\dots,Y_{m})=\nu_c^{\pm}(Y_1,Y_2,\dots,Y_{l},Y_m)+\nu_c^{\pm}(Y_l,\dots,Y_{m}),
\end{aligned}\end{equation}
where  we use that $\nu_c^{\pm}(Y_i,Y_j)=0$ and
\begin{equation}\label{recm3mod}\begin{aligned}
\nu_c^{\pm}(Y_1,Y_2,\dots,Y_m)&=\sum\limits_{j=2}^{m-1}\nu_c^{\pm}(Y_1,Y_j,Y_{j+1})=\sum\limits_{j=1}^{m-2}\nu_c^{\pm}(Y_j,Y_{j+1},Y_m)
\end{aligned}\end{equation}

\par{(ii)} For cyclic sums given by \eqref{cyclsum} we have for  any $2\le l<m$
\begin{equation}\label{rec2}\begin{aligned}
  \mu_c^{\pm}(Y_1,Y_2,\dots,Y_m)&=\mu_c^{\pm}(Y_1,Y_2,\dots,Y_{l})+\mu_c^{\pm}(Y_1,Y_l,\dots,Y_{m})-\rank w(Y_1,Y_l)\\&=\mu_c^{\pm}(Y_1,Y_2,\dots,Y_{l},Y_m)+\mu_c^{\pm}(Y_l,\dots,Y_{m})-\rank w(Y_l,Y_m),
\end{aligned}\end{equation}
where we use that $\mu_c^{\pm}(Y_i,Y_j)=\rank w(Y_i,Y_j)$ and
\begin{equation}\label{recm3}\begin{aligned}
\mu_c^{\pm}(Y_1,Y_2,\dots,Y_m)&=\sum\limits_{j=2}^{m-1}\mu^{\pm}_c(Y_1,Y_j,Y_{j+1})-\sum\limits_{j=3}^{m-1}\rank w(Y_1,Y_j)=\sum\limits_{j=1}^{m-2}\mu_c^{\pm}(Y_j,Y_{j+1},Y_m)-\sum\limits_{j=2}^{m-2}\rank w(Y_j,Y_m)
\end{aligned}\end{equation}
\end{proposition}
\begin{proof} The proof of the equalities in (i) follows from the definition of cyclic sums \eqref{cyclsummod}. Indeed, applying \eqref{cyclsummod} to the cyclic sums $\nu_c^{-}(\cdot)$ in the right-hand sides of \eqref{rec2mod} we see that the addends $-\mu(Y_1,Y_l)$ and $\mu(Y_1,Y_l)$  are cancelled in the first identity, and similarly, $\mu(Y_l,Y_m)$ and $-\mu(Y_l,Y_m)$ are cancelled in the second one. The proof for
$\nu_c^{+}(\cdot)$ is similar.

By a similar way, we prove the first equality in \eqref{recm3mod} cancelling the terms $-\mu(Y_1,Y_{j+1})$ and $\mu(Y_1,Y_{j+1})$ for $j=2,\dots,m-1.$ For the proof of the second one it is sufficient to cancel the terms $-\mu(Y_{j+1},Y_m)$ and $\mu(Y_{j+1},Y_m)$ for $j=1,\dots,m-2.$ Remark  that the second equalities in \eqref{rec2mod} and \eqref{recm3mod} can be also proved via subsequent applications of Proposition~\ref{invar1}(iii) to the first ones.

The proof of (ii) follows from (i) according to connection \eqref{add5} between the cyclic sums  $\nu_c^{\pm}(\cdot)$ and $\mu_c^{\pm}(\cdot)$  in the left and right-hand sides of the equalities in (ii) and (i). The proof is completed.
\end{proof}
\begin{remark}\label{comments}
\par{(i)} The results in Propositions~\ref{invar},~\ref{invar1},~\ref{recrel} were inspired by properties of the Kashiwara index \eqref{kash} (see \cite{Kash,seC.rL.eyM1994}). In particular, by \cite[Property I]{seC.rL.eyM1994} we have $\tau(L_1,L_2,L_3)=\tau(L_3,L_1,L_2)=\tau(L_2,L_3,L_1)$ and
 by \eqref{inv1} we have the similar property for \eqref{cyclsum}
\begin{equation}\label{cyclperm3}\begin{aligned}
  \mu_c^{\pm}(Y_1,Y_2,Y_3)=\mu_c^{\pm}(Y_3,Y_1,Y_2)=\mu_c^{\pm}(Y_2,Y_3,Y_1).
\end{aligned}\end{equation}
At the same time instead of the  property $\tau(L_1,L_2,L_3)=-\tau(L_2,L_1,L_3)$ (see \cite[Property I]{seC.rL.eyM1994}) for $m=3$ we have
\begin{equation}\label{permut2}\begin{aligned}
  \mu_c^{\pm}(Y_1,Y_2,Y_3)=\sum\limits_{i<j}\rank w(Y_i,Y_j)&-\mu_c^{\pm}(Y_2,Y_1,Y_3),\;\mu_c^{\pm}(Y_2,Y_1,Y_3)=\mu_c^{\pm}(Y_3,Y_2,Y_1)=\mu_c^{\pm}(Y_1,Y_3,Y_2),
  \end{aligned}\end{equation}
where we used the connection $\mu_c^{\pm}(Y_2,Y_1,Y_3)=\mu_c^{\pm}(Y_3,Y_2,Y_1)=\mu_c^{\mp}(Y_1,Y_2,Y_3)$ according to Proposition~\ref{invar1}(ii),(iii) and then applied  Proposition~\ref{invar1}(iv).
\par{(ii)} For the cyclic sum $\nu_c^{\pm}(Y_1,Y_2,Y_3)$ one can easily derive from \eqref{permut2}
\begin{equation}\label{permut2nu}\begin{aligned}
  \nu_c^{\pm}(Y_1,Y_2,Y_3)&=\rank w(Y_1,Y_2)-\nu_c^{\pm}(Y_2,Y_1,Y_3)=\rank w(Y_2,Y_3)-\nu_c^{\pm}(Y_1,Y_3,Y_2)\\&=\rank w(Y_1,Y_2)+\rank w(Y_2,Y_3)-\rank w(Y_1,Y_3)-\nu_c^{\pm}(Y_3,Y_2,Y_1),
  \end{aligned}\end{equation}
where we use \eqref{add5} incorporating the order of the components in the cyclic sums.
\par{(iii)} Observe that the cyclic sums \eqref{cyclsummod} do not obey  property \eqref{inv2},  but one can easily derive the connections
\begin{equation*}\begin{aligned}
  \nu_c^{\pm}(Y_1,Y_2,\dots,Y_m)&=\mu_c^{\pm}(Y_1,Y_2,\dots,Y_m)-\rank w(Y_m,Y_1)=\mu_c^{\pm}(Y_m,Y_1,\dots,Y_{m-1})-\rank w(Y_m,Y_1)\\&=\nu_c^{\pm}(Y_m,Y_1,\dots,Y_{m-1})-\rank w(Y_m,Y_1)+\rank w(Y_{m-1},Y_m)\\&=\nu_c^{\pm}(Y_{m-1},Y_m,Y_1,\dots,Y_{m-2})-\rank w(Y_m,Y_1)+\rank w(Y_{m-2},Y_{m-1})\\&\quad \dots \quad \dots\\&=\nu_c^{\pm}(Y_2,Y_3,\dots,Y_m,Y_1)-\rank w(Y_m,Y_1)+\rank w(Y_{1},Y_{2})
\end{aligned}\end{equation*}
where we used  \eqref{add5} and applied Proposition~\ref{invar1}(iii).
\par{(iv)} Relations \eqref{recm3mod} for $m=4$ imply the following analogs of the "cocycle condition", see \cite[Theorem A.3.2(ii)]{Kash}
\begin{equation}\label{Horprop1}
  \nu^{\pm}_c(Y_1,Y_2,Y_3,Y_4)=\nu^{\pm}_c(Y_1,Y_2,Y_3)+\nu^{\pm}_c(Y_1,Y_3,Y_4)=\nu^{\pm}_c(Y_1,Y_2,Y_4)+\nu^{\pm}_c(Y_2,Y_3,Y_4),
\end{equation}
where \eqref{Horprop1} can be rewritten in the form
\begin{equation}\label{Horprop}\begin{aligned}
  \nu^{\pm}_c(Y_1,Y_2,Y_4)-\nu^{\pm}_c(Y_1,Y_2,Y_3)=\nu^{\pm}_c(Y_1,Y_3,Y_4)-\nu^{\pm}_c(Y_2,Y_3,Y_4).
\end{aligned}\end{equation}
\end{remark}
\section{Index results for cyclic sums}\label{sec2}
In this section we present the proof of Theorem~\ref{maincon} which connects $\mu_c^{\pm}(Y_1,Y_2,\dots,Y_m)$ with the negative and positive inertia $i_{\mp}$ of the $(m n)\times (mn)$ symmetric matrix given by \eqref{Sm}. Based on this result we prove a similar connection for $\nu_c^{\pm}(Y_1,Y_2,\dots,Y_m).$

The proof of Theorem~\ref{maincon} is based on the fact that for  Lagrangian subspaces $L_1,\,L_2,\,\dots,\,L_m$ there exists a Lagrangian subspace $L_{R}$ such that $L_j\cap L_R=\{0\},$ see e.g. \cite{Kash}, \cite{seC.rL.eyM1994}. For completeness we prove a similar result for the frames $Y_1,\,Y_2,\dots,Y_m$ presenting  a special transformation matrix associated with $L_R$ (see also \cite{Abramov2001}, where this result is used  for $m=1$ and \cite[Lemma 3.5]{jE20} for the case $m=2$).
\begin{lemma}\label{simultan}
Let $Y_k=\binom{X_k}{U_k},\,k=1,2,\dots,m$ be $2n\times n$ matrices with condition \eqref{conj}. Then there exists $\alpha\in\mathbb{R},\,\alpha\ne \pi k/2,\;k\in \mathbb{Z}$ such that for the transformation matrix
\begin{equation}\label{angle}
  R_{\alpha}=\begin{pmatrix}
                       \cos (\alpha)\,I & \sin (\alpha)\,I \\
                       -\sin (\alpha)\,I & \cos (\alpha)\,I \\
                     \end{pmatrix}
 \end{equation}
 we have
 \begin{equation}\label{nonsing}
   \det \tilde X_k \ne 0,\quad k=1,2,\dots,m,\;\tilde X_k=\cos(\alpha)X_k-\sin(\alpha)U_k
 \end{equation}
 where $\tilde X_k,\,k=1,2,\dots,m$ are the upper blocks of $\tilde Y_k=R_{\alpha}^{-1}Y_k.$
\end{lemma}
\begin{proof}
Consider the determinants
\begin{equation*}
 P_k(\gamma)= \det(X_k-\gamma U_k),\,k=1,2,\dots,m
\end{equation*}
for a complex parameter $\gamma\in\mathbb{C}.$ All these determinants have a polynomial dependence in  $\gamma\in\mathbb{C}.$ Since these functions are nontrivial (for example, for $\gamma:=i$) conditions \eqref{nonsing} are satisfied for all  $\gamma_0\in\bR$ which do not coincide with real roots of the polynomials  $\det(X_k-\gamma U_k),\,k=1,2,\dots,m$. Remark that for $m<\infty$ we have a finite number of the roots of $\det(X_k-\gamma U_k),\,k=1,2,\dots,m$, then such $\gamma_0\in\bR$ does exists.  Remark also that one can put $\alpha:=\arctan(\gamma_0)+\pi k,$ for $R_{\alpha}$  given by \eqref{angle}. In this case condition $\alpha\ne \pi k,\,k\in \mathbb{Z}$ is satisfied for $\gamma_0\ne 0.$ The proof is completed.
\end{proof}
Applying Lemma~\ref{simultan} and Proposition~\ref{invar} we derive the following representations for \eqref{cyclsum} and \eqref{cyclsummod}
\begin{lemma}\label{main}
Let $R_{\alpha}$ be chosen according to \eqref{nonsing}, then for the cyclic sums defined by \eqref{cyclsum} and \eqref{cyclsummod} we have
\begin{equation}\label{durtrans}
  \mu_c^{\pm}(Y_1,Y_2,\dots,Y_m)=\mu_c^{\pm}(R_{\alpha}^{-1}Y_1,R_{\alpha}^{-1}Y_2,\dots,R_{\alpha}^{-1}Y_m)=\sum\limits_{k=1}^{m-1}\ind(\pm(\tilde Q_{k}-\tilde Q_{k+1}))+\ind(\pm(\tilde Q_m-\tilde Q_1))
\end{equation}
and similarly
\begin{equation}\label{durtrans*}
  \nu^{\pm}_c(Y_1,Y_2,\dots,Y_m)=\nu^{\pm}_c(R_{\alpha}^{-1}Y_1,R_{\alpha}^{-1}Y_2,\dots,R_{\alpha}^{-1}Y_m)=\sum\limits_{k=1}^{m-1}\ind(\pm(\tilde Q_{k}-\tilde Q_{k+1}))-\ind(\pm(\tilde Q_1-\tilde Q_m)),
\end{equation}
where
\begin{equation}\label{trQ}
  \tilde Q_k=\tilde Q_k^T=\tilde U_k \tilde X_k^{-1},\,\tilde Y_k=R_{\alpha}^{-1}Y_k=\binom{\tilde X_k}{\tilde U_k}.
\end{equation}
\end{lemma}
\begin{proof}
The proof follows from Proposition~\ref{invar}, \eqref{cyclsum}, \eqref{cyclsummod}, and the definition of the comparative indices $$\mu(\tilde Y_k,\tilde Y_{k+1})=\ind(\tilde Q_{k+1}-\tilde Q_k),\,k=1,\dots,m,\; \mu(\tilde Y_m,\tilde Y_1)=\ind(\tilde Q_{1}-\tilde Q_m)$$ and the dual comparative indices $\mu^*(\tilde Y_k,Y_{k+1})=\ind(\tilde Q_{k}-\tilde Q_{k+1}),\,\mu^*(\tilde Y_m,\tilde Y_1)=\ind(\tilde Q_{m}-\tilde Q_1)$ under the nonsingularity condition \eqref{nonsing} according to \eqref{compnonsing}.
\end{proof}

Now  we present the proof of Theorem~\ref{maincon}.

\begin{proof}[Proof of Theorem~\ref{maincon}]
Consider  the matrix $S_{1,2,\dots,m}$ given by \eqref{Sm}. By Lemmas~\ref{simultan},~\ref{main}  we have for the Wronskians $w(Y_i,Y_j)$ in \eqref{Sm}
 \begin{equation}\label{repWR2}
   w(Y_i,Y_j)=w(R_{\alpha}Y_i,R_{\alpha}Y_j)=w(\tilde Y_i,\tilde Y_j)=\tilde X_i^T(\tilde Q_j-\tilde Q_i)\tilde X_j,\;\tilde Q_k=\tilde U_k \tilde X_k^{-1},\,k=i,j.
 \end{equation}
Hence $S_{1,2,\dots,m}$ can be rewritten in the form
\begin{equation}\label{congr}
  S_{1,2,\dots,m}=\diag\{\tX_1^T,\;\tX_2^T,\dots,\tX_m^T\}\tilde S_{1,2,\dots,m}\diag\{\tX_1,\;\tX_2,\dots,\tX_m\},
\end{equation}
where $\tilde S_{1,2,\dots,m}$ consists of $n\times n$ blocks $ \tilde S_{1,2,\dots,m}(i,j),\,\;i,j=1,2,\dots,m$
\begin{equation}\label{SQm}
  \tilde S_{1,2,\dots,m}(i,j)=\tQ_j-\tQ_i,\quad \tilde S_{1,2,\dots,m}(j,i)=\tilde S_{1,2,\dots,m}^T(i,j),\;j\ge i.
\end{equation}
Here the symmetric matrices $\tQ_k$ are given by \eqref{trQ} and the matrices $\tX_j,\,j=1,2,\dots,m$ are nonsingular according to \eqref{nonsing} in Lemma~\ref{simultan}. Then we have
\begin{equation*}
  \ind(\pm S_{1,\dots,m})=\ind(\pm \tilde S_{1,\dots,m}).
\end{equation*}

Introduce the matrices $M_m,\,m\ge 2$ with $n\times n$ blocks such that $M_{i,i}=I,\,\;M_{i,i+1}=-I,\;M_{m,1}=I,$ and $M_{i,j}=0$ otherwise.
In particular, we have
\begin{equation*}\label{exM234}
  M_2=\begin{pmatrix}
        I & -I \\
        I & I \\
      \end{pmatrix},\;M_3=\begin{pmatrix}
                            I & -I & 0 \\
                            0 & I & -I \\
                            I & 0 & I \\
                          \end{pmatrix},\,M_4=\begin{pmatrix}
                                                I & -I & 0 & 0 \\
                                                0 & I & -I & 0 \\
                                                0& 0& I& -I \\
                                                I & 0 & 0 & I \\
                                              \end{pmatrix}.
\end{equation*}
We prove that  the matrices $M_m$ are nonsingular using their partitioned form
\begin{equation*}\label{M}
  M_m=\begin{pmatrix}
                             L & N \\
                             K & I \\
                           \end{pmatrix},\,L=\begin{pmatrix}
                                                     I & -I & 0 & \dots & \dots & 0 \\
                                                     0 & I & -I & 0 & \dots & 0 \\
                                                     \dots & \dots & \dots & \dots & \dots &\dots \\
                                                     0 & 0 & \dots & 0 & I &-I \\
                                                     0 & 0 & \dots & 0 & 0& I \\
                                                   \end{pmatrix},\quad N=(0\,0\,\dots\,0\,-I)^T,\;K=(I\,0\,\dots\,0).
  \end{equation*}
Then, by \cite[Proposition 2.8.3]{Bernstein} $\det M_m=\det L\det(I-K L^{-1} N)=\det(2I)=2^n,$ where we
used that $L^{-1}$ has the $n\times n$ blocks $L^{-1}(i,j)=I,\,j\ge i$ and $L^{-1}(i,j)=0$ otherwise.

For arbitrary $m\ge 2$ we have for $\tilde S_{1,\dots,m}$ given by \eqref{congr}, \eqref{SQm}
\begin{equation}\label{matrBmd}
M_{m} \tilde S_{1,\dots,m} M^{T}_{m}=\hat S_{1,2,\dots,m} ,\quad \hat S_{1,2,\dots,m}=\diag\{2(\tQ_1-\tQ_2),2(\tQ_2-\tQ_3),\dots,2(\tQ_m-\tQ_1)\}.
\end{equation}
Indeed, using the simple structure of $M_m$ we have for the $n\times n$ blocks of  $M_m  \tilde S_{1,\dots,m}$
\begin{equation*}
  M_m  \tilde S_{1,\dots,m}(k,j)=\tQ_k-\tQ_{k+1},\,j\le k,\quad M_m  \tilde S_{1,\dots,m}(k,j)=\tQ_{k+1}-\tQ_k,\,j> k,\quad  k=1,\dots,m-1,\;j=1,\dots,m,
\end{equation*}
and $M_m  \tilde S_{1,\dots,m}(m,j)=\tQ_m-\tQ_1,\,j=1,\dots,m.$ Multiplying $M_m  \tilde S_{1,\dots,m}$ and $M_m^T$ we derive \eqref{matrBmd}.

Applying Lemma~\ref{main} we have by \eqref{matrBmd}
\begin{equation*}\begin{aligned}
  \mu_c^{\pm}(Y_1,Y_2,\dots,Y_m)&=\sum\limits_{k=1}^{m-1}(\ind(\pm(\tilde Q_{k}-\tilde Q_{k+1}))+\ind(\pm(\tilde Q_m-\tilde Q_1))=\ind(\pm \hat S_{1,2,\dots,m})= \ind(\pm \tilde S_{1,2,\dots,m})=\ind(\pm S_{1,2,\dots,m}).
\end{aligned}\end{equation*}
The proof is completed.
\end{proof}

Introduce the notation for the blocks of $S_{m,1,2,\dots,m-1}$ defined by \eqref{Sm} for the cyclic permutation $Y_m,\,Y_1,\,Y_2,\dots,\,Y_{m-1}$
\begin{equation}\label{imprnot}\begin{aligned}
 S_{m,1,2,\dots,m-1}&=\begin{pmatrix}
                                          0 &  \mathcal W \\
                                          \mathcal  W^{T} & S_{1,2,\dots,m-1} \\
                                        \end{pmatrix},\,m\ge 2,\;\mathcal W =(w_{m,1}\; w_{m,2}\;\dots\;w_{m,m-1}),\\S_{m,1,2,\dots,m-1}&=\begin{pmatrix}
                       0 & w_{m,1} & \mathcal N\\
                       w_{m,1}^T & 0 & \mathcal K \\
                       \mathcal N^T & \mathcal K^T & S_{2,3,\dots,m-1} \\
                     \end{pmatrix},\,m\ge 3,  \;\mathcal N=(w_{m,2}\;w_{m,3}\,\dots\;w_{m,m-1}),\;\mathcal K= (w_{1,2}\; w_{1,3}\;\dots\;w_{1,m-1}),
\end{aligned}\end{equation}
where in the first representation $\mathcal W \in \mathbb{R}^{n\times(m-1)n},\,m\ge 2$ and for $m=2$ we put $S_1:=0_n.$ We also have $\mathcal K,\,\mathcal N \in \mathbb{R}^{n\times(m-2)n},\,m\ge 3,$ and  for $m=3$ we put $S_2:=0_n.$

Applying  \cite[Theorem 2.3]{Tian}, see \eqref{tianin}, we derive the following corollary to Theorem~\ref{maincon}.
\begin{corollary}\label{maincon1+}
Under notation \eqref{imprnot} we have the following representations for cyclic sums \eqref{cyclsum}
\begin{equation}\label{maincon1a}\begin{aligned}
  \mu_c^{\pm}(Y_1,Y_2,\dots,Y_m)&=\rank \mathcal W+\ind(\pm F_{\mathcal W}\;S_{1,2,\dots,m-1}\;F_{\mathcal W}),\quad m\ge 2,\end{aligned}\end{equation}

  \begin{equation}\label{maincon1b}\begin{aligned}
  \mu_c^{\pm}(Y_1,Y_2,\dots,Y_m)&=\rank w_{m,1}+\ind(\pm \bar S_{1,2,\dots,m-1}),\quad \bar S_{1,2,\dots,m-1}= \begin{pmatrix}
                       0 &\tilde {\mathcal M} \\
                      \tilde{\mathcal M}^T  & S_{2,3,\dots,m-1}- \mathcal D-\mathcal D^T\\[2mm]
                     \end{pmatrix},\,m\ge 3,\\ \tilde{\mathcal M}&=E_{w_{m,1}}\mathcal N=(I-w_{m,1}w_{m,1}^{\dag}) (w_{m,2}\,w_{m,3}\,\dots \,w_{m,m-1}),\quad \mathcal D=\mathcal K^T \,w_{m,1}^{\dag} \mathcal N,\end{aligned}\end{equation}
where the matrix $\tilde{\mathcal M}$ can be replaced by
\begin{equation}\label{mtgen}
 {\mathcal M}=F_{w_{m,1}}\mathcal K=(I-w_{m,1}^{\dag}w_{m,1}) (w_{1,2}\,w_{1,3}\,\dots \,w_{1,m-1}).
\end{equation}
\end{corollary}
\begin{proof}
Applying Proposition~\ref{invar1}(ii) we have by Theorem~\ref{maincon}
\begin{equation}\label{permmaincon}
 \mu_c^{\pm}(Y_1,Y_2,\dots,Y_m)=\mu_c^{\pm}(Y_m,Y_1,\dots,Y_{m-1}) =\ind(\pm S_{m,1,2,\dots,m-1}),
\end{equation}
where $S_{m,1,2,\dots,m-1}$ is given by the first matrix in \eqref{imprnot}. Then  \eqref{maincon1a} is derived by the direct application of the second formula in \eqref{tianin} to \eqref{permmaincon}.

Consider the proof of  \eqref{maincon1b}. Rewrite $S_{m,1,2,\dots,m-1}$ in the second equality \eqref{imprnot} as follows
\begin{equation}\label{add9}
S_{m,1,2,\dots,m-1}=\begin{pmatrix}
                      A & B \\
                      B^T & D \\
                    \end{pmatrix},\;
  A:=\begin{pmatrix}
                    0 &  w_{m,1} \\
                    w_{m,1}^T & 0 \\
                  \end{pmatrix},\quad B:=\begin{pmatrix}
                                       \mathcal N \\
                                        \mathcal K\\
                                     \end{pmatrix},\quad
                  D:=S_{2,3,\dots,m-1}
\end{equation}
and then apply the first formula in \eqref{tianin} incorporating that
                    $\ind(\pm A)=\rank w_{m,1},$
see  \eqref{n=2}. We also compute by \eqref{tianin}
\begin{equation*}
  D-B^TA^{\dag}B=S_{2,3,\dots,m-1}-\mathcal K^T \,w_{m,1}^{\dag} \mathcal N-\mathcal N^T \,w_{m,1}^{\dag\;T} \mathcal K=S_{2,3,\dots,m-1}- \mathcal D-\mathcal D^T
\end{equation*}
and
\begin{equation}\label{add10}
  E_A B=\diag\{E_{w_{m,1}},\,E_{w_{m,1}^T}\}\begin{pmatrix}
                                       \mathcal N \\
                                        \mathcal K\\
                                     \end{pmatrix}=
  \begin{pmatrix}
          \tilde{\mathcal M} \\
          {\mathcal M} \\
        \end{pmatrix}.
\end{equation}
A key step in the proof is connected with  properties of $\tilde{\mathcal M}$ and ${\mathcal M}$ defined in \eqref{maincon1b}, \eqref{mtgen}.
We have that there exists a  nonsingular matrix $\mathcal L$ such that
\begin{equation}\label{propM}
  \mathcal M=\mathcal L\, \tilde{\mathcal M}.
\end{equation}
This fact follows directly from the properties of the comparative index, see the proof of Theorem~2.1 in \cite{jE07}, the proof of Theorem~3.2 in \cite[p. 151]{DEH}. Indeed, it is sufficient to consider the comparative indices $\mu(Z_m^{-1}Y_1, Z_m^{-1}Y_j)$ for $j=2,3,\dots,m-1$ and use the definition of the matrices $\mathcal M_j$ and  $\tilde {\mathcal M}_j$ according to \eqref{compind1} and \eqref{compind} incorporating \eqref{iverseandwrons}
\begin{equation*}
  \mathcal M_{j}=(I-w_{m,1}^{\dag}w_{m,1})w_{1,j},\quad \tilde {\mathcal M}_{j}=(I-w_{m,1}w_{m,1}^{\dag})w_{m,j}.
\end{equation*}
Applying the result from \cite[p. 151]{DEH} we see that $\mathcal M_{j}=\mathcal L \,\tilde {\mathcal M}_{j},$ where $\mathcal L$ depends only on the blocks of $Z_m^{-1}Y_1,$ see \cite[formula (1.74)]{DEH}. Since
\begin{equation*}
  {\mathcal M}=({\mathcal M}_{2}\,{\mathcal M}_{3}\,\dots \,{\mathcal M}_{m-1}),\quad \tilde {\mathcal M}=(\tilde{\mathcal M}_{2}\,\tilde{\mathcal M}_{3}\,\dots \,\tilde{\mathcal M}_{m-1})
\end{equation*}
 formula \eqref{propM} is proved.

By \eqref{propM} the matrix $\begin{pmatrix}
                               0 & E_A B \\
                               (E_A B)^T & D-B^TA^{\dag}B \\
                             \end{pmatrix}$
with  blocks \eqref{add9}, \eqref{add10} in the right-hand side of the first equality in \eqref{tianin}  can be simplified by deleting rows and columns  containing  ${\mathcal M},\,{\mathcal M}^T$ (or $\tilde {\mathcal M},\,\tilde{\mathcal M}^T$). So we see that the dimension of $\bar S_{1,2,\dots,m-1}$ in \eqref{maincon1b} is equal to $(m-1)n\times (m-1)n.$

Finally, we also showed  by \eqref{propM} that $\tilde{\mathcal M}$ can be replaced by ${\mathcal M}.$ The proof is completed.
\end{proof}
Based on Corollary~\ref{maincon1+} we derive the following representations for the cyclic sums $\nu_c^{\pm}(Y_1,Y_2,\dots,Y_m).$
\begin{theorem}\label{maincon23}
Under the notation of Corollary~\ref{maincon1+} we have for  \eqref{cyclsummod}
\begin{equation}\label{maincon2}\begin{aligned}
                     \nu_c^{\pm}(Y_1,Y_2,\dots,Y_m)&=\rank (\tilde{\mathcal M}) +\ind(\pm F_{\mathcal W}\;S_{1,2,\dots,m-1}\;F_{\mathcal W}),\,m\ge 2,\\[2mm]
  \nu_c^{\pm}(Y_1,Y_2,\dots,Y_m)&=\ind(\pm \bar S_{1,2,\dots,m-1})=\rank (\tilde{\mathcal M})+\ind(\pm F_{\tilde{\mathcal M}} (S_{2,3,\dots,m-1}-\mathcal D-\mathcal D^T )F_{\tilde{\mathcal M}}) ,\;m\ge 3,
\end{aligned}\end{equation}
where $\bar S_{1,2,\dots,m-1}$, $\tilde{\mathcal M}$ are given by \eqref{maincon1b} and $\tilde{\mathcal M}$ can be replaced by $ {\mathcal M}$ defined by \eqref{mtgen}.
\end{theorem}
\begin{proof}
Representation  \eqref{maincon1a}  implies the first connection in \eqref{maincon2}, where we use \eqref{add5} and apply \eqref{stian1}
\begin{equation}\label{maincon5}
  \rank \mathcal W=\rank w_{m,1}+\rank (\tilde{\mathcal M}).
\end{equation}

By a similar way, the equality $\nu_c^{\pm}(Y_1,Y_2,\dots,Y_m)=\ind(\pm \bar S_{1,2,\dots,m-1})$  in \eqref{maincon2} follows from  \eqref{maincon1b}, where we use  connection \eqref{add5} between $\nu_c^{\pm}(Y_1,Y_2,\dots,Y_m)$ and $\mu_c^{\pm}(Y_1,Y_2,\dots,Y_m).$ Next we compute the index $\ind(\pm \bar S_{1,2,\dots,m-1})$ applying the second identity in \eqref{tianin}. Observe that the matrix $F_{\tilde{\mathcal M}}$ stays the same after the replacement $\tilde{\mathcal M}$ by ${\mathcal M}$ because of connection \eqref{propM}, see \cite[Theorem 8, Lemma 3]{BenIasr}. The proof is completed.
\end{proof}
\begin{remark}\label{rem1}
\par{(i)} It follows from \eqref{maincon2} that for $m\ge 3$
\begin{equation}\label{maincon3}\begin{aligned}
  \ind(\pm F_{\mathcal W}\;S_{1,2,\dots,m-1}\;F_{\mathcal W})=\ind(\pm F_{\tilde{\mathcal M}} (S_{2,3,\dots,m-1}-\mathcal D-\mathcal D^T)F_{\tilde{\mathcal M}}),
\end{aligned}\end{equation}
in particular, for $m=3$ formula \eqref{maincon3} presents the second addend in the definition of the comparative index $\mu(Z_3^{-1}Y_1, Z_3^{-1}Y_2),$ compare with \eqref{compwron2}.
\par{(ii)} Summing $\mu_c^{-}(\cdot)$ and $\mu_c^{+}(\cdot)$ given by \eqref{maincon1+} and incorporating \eqref{maincon3}  we have for $m\ge 3$
\begin{equation}\label{maincon4}\begin{aligned}
  \rank(\pm F_{\mathcal W}\;S_{1,2,\dots,m-1}\;F_{\mathcal W})&=\rank(\pm F_{\tilde{\mathcal M}} (S_{2,3,\dots,m-1}-\mathcal D-\mathcal D^T)F_{\tilde{\mathcal M}})\\=\sum\limits_{j=1}^{m-1}\rank w(Y_j,Y_{j+1})+\rank w(Y_m,Y_1)&-2\rank \mathcal W ,\quad \mathcal W =(w_{m,1}\; w_{m,2}\;\dots\;w_{m,m-1}),
\end{aligned}\end{equation}
where we also used Proposition~\ref{invar1}(iv). Observe that  $\rank \mathcal W$ can be presented in terms of the dimension of $L_1+L_2+\dots+L_m$ (or $L_1\cap L_2\cap\dots\cap L_m$)
\begin{equation}\label{dimrank}\begin{aligned}
  \dim(L_1+L_2+\dots+L_m)&=\rank (Y_1\,Y_2\,\dots Y_m)=\rank Z_m(Z_m^{-1}Y_1\,Z_m^{-1}Y_2\,\dots\,Z_m^{-1}Y_{m-1}\,\;\bnn)\\&=\rank (Z_m^{-1}Y_1\,Z_m^{-1}Y_2\,\dots\,Z_m^{-1}Y_{m-1}\,\bnn)=n+\rank \mathcal W,
\end{aligned}\end{equation}
where we used \eqref{matr} and \eqref{iverseandwrons}. Finally, substituting $\rank \mathcal W=\dim(L_1+L_2+\dots+L_m)-n$ from \eqref{dimrank} and $\rank w_{i,j}=n-\dim(L_i\cap L_j)$ into \eqref{maincon4} we derive by $\dim(L_1\cap L_2\cap\dots \cap\,L_m)=2n-\dim(L_1+L_2+\dots+L_m)$
\begin{equation}\label{geom}\begin{aligned}
  \rank(\pm F_{\mathcal W}\;S_{1,2,\dots,m-1}\;F_{\mathcal W})&=\rank(\pm F_{\tilde{\mathcal M}} (S_{2,3,\dots,m-1}-\mathcal D-\mathcal D^T)F_{\tilde{\mathcal M}})\\=2\dim(L_1\cap L_2\cap\dots \cap\,L_m)+(m-2)n&-\sum\limits_{j=1}^{m-1}\dim(L_j\cap L_{j+1})-\dim (L_m\cap L_1)\ge 0,
\end{aligned}\end{equation}
compare with \cite[Corollary 3.7]{Zhou} for the case $m=3.$
\par{(iii)} One can verify that  diagonal blocks of $\mathcal D$ defined in \eqref{maincon1b} are symmetric on the image of $F_{\tilde{\mathcal M}}.$ We incorporated this fact deriving the right-hand side of \eqref{maincon3} for $m=3,$ see Example~\ref{ex1}. \end{remark}
From Theorem~\ref{maincon} and Lemma~\ref{prop ind}(vi) we also derive the following estimates for cyclic sums \eqref{cyclsum}, \eqref{cyclsummod}.
\begin{corollary}\label{estimcycls}
For the cyclic sums \eqref{cyclsum}, \eqref{cyclsummod} we have the estimates
\begin{equation}\label{estcycl}\begin{aligned}
0& \le r\le
\mu_c^{\pm}(Y_1,Y_2,\dots,Y_m)\le P,\quad 0 \le r-\rank w_{1,m}\le
\nu_c^{\pm}(Y_1,Y_2,\dots,Y_m)\le P-\rank w_{1,m},\\r:&=\max\limits_{l<k,\;l,k=1,\dots,m}(\rank w(Y_l,Y_{k})),\\P:&=\sum\limits_{j=1}^{m-1}\min(\rank w(Y_j,Y_{j+1}),\rank w(R\bnn,Y_{j+1}))+\min(\rank w(Y_m,Y_{1}),\rank w(R\bnn,Y_{1})),
\end{aligned}\end{equation}
where $R$ is  arbitrary symplectic matrix.
\end{corollary}
\begin{proof}
 By Proposition~\ref{invar} we have $\mu_c^{\pm}(Y_1,Y_2,\dots,Y_m)=\mu_c^{\pm}(R^{-1}Y_1,R^{-1}Y_2,\dots,R^{-1}Y_m)$, then the upper bounds   in \eqref{estcycl} follow from Lemma~\ref{prop ind}(vi) applied to the sums of the comparative indices $\mu(R^{-1}Y_j,R^{-1}Y_{j+1}),\,j=1,\dots,m-1,$ $\mu(R^{-1}Y_m,R^{-1}Y_1)$ according to  Definition~\ref{cyclsumd} and \eqref{add5}.

For the proof of the lower bounds we apply Theorem~\ref{maincon}. By \eqref{Sm} one can chose the $2n\times 2n$ principal  submatrix of $S_{1,2,\dots,m}$ depending on the Wronskians $w(Y_i,Y_j),\,i<j$ in the form $S_{i,j}=\smat{0}{w_{i,j}}{w_{i,j}^T}{0}$
with $\ind S_{i,j}=\ind (-S_{i,j})=\rank w(Y_i,Y_j),$ (see \cite[Fact 5.8.8]{Bernstein}), then the lower bound in the first inequality \eqref{estcycl} follows from the inequalities
\begin{equation*}
 \rank w(Y_i,Y_j)  =\ind (\pm S_{i,j})\le \ind (\pm S_{1,2,\dots,m})=\mu^{\pm}_c(Y_1,Y_2,\dots,Y_m),
\end{equation*}
see  \cite[Fact 5.8.20]{Bernstein}. The lower bound for \eqref{cyclsummod} follows from \eqref{add5}.
\end{proof}
\section{Applications}\label{sec3}
\subsection{Connections with the Kashiwara index}
Recall that according to \cite{Kash}  (see Section~\ref{sec0})  {the Kashiwara index} $\tau(L_1,L_2,L_3)$ is defined as the signature  $ \tau(L_1,L_2,L_3)=\sign(B(x,x))=i_{+}(B(x,x))-i_{-}(B(x,x))$
of the quadratic form $B(x,x)$ given by \eqref{bilin}
 \begin{equation*}
 B(x,x):=B((x_1,x_2,x_3), (x_1,x_2,x_3))=w(x_1,x_2)+w(x_2,x_3)+w(x_3,x_1)
\end{equation*}
defined on $(x_1,x_2,x_3)\in L_1\oplus L_2\oplus L_3.$  This definition is generalized to the case $m\ge 3$  according to \eqref{kashm}.

The main result of this section is the following theorem.
\begin{theorem}\label{kashcycl}
We have the following connections for $m=3$
\begin{equation}\label{kashcycls3+mod}
 \tau(L_1,L_2,L_3)= \mu_c^{+}(Y_1,Y_2,Y_3)-\mu_c^{-}(Y_1,Y_2,Y_3)=\nu_c^{+}(Y_1,Y_2,Y_3)-\nu_c^{-}(Y_1,Y_2,Y_3).
\end{equation}

For the case $m\ge 3$ we have  for \eqref{cyclsum}, \eqref{cyclsummod}, and \eqref{kashm}
\begin{equation}\label{kashcyclsm+m}
 \tau(L_1,L_2,\dots,L_m)= \mu_c^{+}(Y_1,Y_2,\dots,Y_m)-\mu_c^{-}(Y_1,Y_2,\dots,Y_m)=\nu_c^{+}(Y_1,Y_2,\dots,Y_m)-\nu_c^{-}(Y_1,Y_2,\dots,Y_m),
\end{equation}
\begin{equation}\label{muthrkash}\begin{aligned}
  \mu_c^{\pm}(Y_1,Y_2,\dots,Y_m)&=\frac{1}{2}(\sum\limits_{j=1}^{m-1}\rank w(Y_j,Y_{j+1})+\rank w(Y_m,Y_1)\pm \tau(L_1,L_2,\dots,L_m)),\\\nu_c^{\pm}(Y_1,Y_2,\dots,Y_m)&=\frac{1}{2}(\sum\limits_{j=1}^{m-1}\rank w(Y_j,Y_{j+1})-\rank w(Y_m,Y_1)\pm \tau(L_1,L_2,\dots,L_m)).
\end{aligned}\end{equation}
\end{theorem}
\begin{proof}
  Introduce the matrix of the quadratic form $B(x,x)$
\begin{equation}\label{matrB3}\begin{aligned}
  S_{B}&=\frac{1}{2}\begin{pmatrix}
        0 & w(Y_1,Y_2) & -w(Y_1,Y_3) \\
        -w(Y_2,Y_1) & 0 & w(Y_2,Y_3) \\
        w(Y_3,Y_1) & -w(Y_3,Y_2) & 0 \\
      \end{pmatrix}=\frac{1}{2}K\begin{pmatrix}
        0 & -w(Y_1,Y_2) & -w(Y_1,Y_3) \\
      - w^T(Y_1,Y_2) & 0 & -w(Y_2,Y_3) \\
       -w^T(Y_1,Y_3)& -w^T(Y_2,Y_3) & 0 \\
      \end{pmatrix}K\\&=-\frac{1}{2}KS_{123}K,\quad K=\diag\{I,-I,I\},
\end{aligned}\end{equation}
where the matrix $S_{123}$ is given by \eqref{Sm} for $m=3$. Then we have
\begin{equation*}\begin{aligned}
\tau(L_1,L_2,L_3)=\sign(B(x,x))=i_{+}(B(x,x))-i_{-}(B(x,x))=\ind(-S_B)-\ind(S_B)=\ind(S_{123})-\ind(-S_{123})
\end{aligned}\end{equation*}
and by Theorem~\ref{maincon} for $m=3$ we derive the first equality in \eqref{kashcycls3+mod}. The second one follows from \eqref{add5} for $m=3.$

For the proof of \eqref{kashcyclsm+m} we use Proposition~\ref{recrel}. By \eqref{add5},  \eqref{recm3mod} we have
\begin{equation*}\begin{aligned}
  \mu_c^{+}(Y_1,Y_2,\dots,Y_m)&-\mu_c^{-}(Y_1,Y_2,\dots,Y_m)=\nu_c^{+}(Y_1,Y_2,\dots,Y_m)-\nu_c^{-}(Y_1,Y_2,\dots,Y_m)\\
  &=\sum\limits_{j=2}^{m-1}(\nu_c^{+}(Y_1,Y_j,Y_{j+1})-\nu_c^{-}(Y_1,Y_j,Y_{j+1}))=\sum\limits_{j=2}^{m-1}\tau(L_1,L_j,L_{j+1})=\tau(L_1,L_2,\dots,L_m),
\end{aligned}\end{equation*}
where we also used \eqref{kashcycls3+mod} for the case $\nu_c^{+}(Y_1,Y_j,Y_{j+1})-\nu_c^{-}(Y_1,Y_j,Y_{j+1})=\tau(L_1,L_j,L_{j+1}).$ The proof of \eqref{kashcyclsm+m} is completed.

Finally, by Proposition~\ref{invar1}(iv) (see \eqref{con1+}) we have
\begin{equation*}\begin{aligned}
  \mu_c^{+}(Y_1,Y_2,\dots,Y_m)+\mu_c^{-}(Y_1,Y_2,\dots,Y_m)=\sum\limits_{j=1}^{m-1}\rank w(Y_j,Y_{j+1})+\rank w(Y_m,Y_1),
\end{aligned}\end{equation*}
then by summing (subtracting) the last equality and $\mu_c^{+}(Y_1,Y_2,\dots,Y_m)-\mu_c^{-}(Y_1,Y_2,\dots,Y_m)=\tau(L_1,L_2,\dots,L_m)$ we derive the first equality in \eqref{muthrkash}. The second one then follows from \eqref{add5}. The proof is completed.
\end{proof}
Based on connections \eqref{muthrkash} some properties of \eqref{cyclsum}, \eqref{cyclsummod} can be derived from similar properties of the Kashiwara index. So we have the following corollary to Theorem~\ref{kashcycl}.
\begin{corollary}\label{dependrank}
 Suppose that $Y_1(t),\,Y_2(t),\dots,Y_m(t)$ with conditions \eqref{conj}  are continuous functions of $t\in[a,b].$ Then under the assumption
\begin{equation}\label{constdim}
  \rank w(Y_j(t),Y_{j+1}(t))=\const,\,j=1,2,\dots,m-1,\quad \rank w(Y_1(t),Y_{m}(t))=\const,\,t\in[a,b]
\end{equation}
we have
\begin{equation}\label{constvar}
  \mu_c^{\pm}(Y_1(t),\,Y_2(t),\dots,Y_m(t))=\const,\quad \nu_c^{\pm}(Y_1(t),\,Y_2(t),\dots,Y_m(t))=\const.
\end{equation}
\end{corollary}
\begin{proof}
Under assumption \eqref{constdim} the Kashiwara index $\tau(L_1(t),L_2(t),\dots,L_m(t))$ (where the Lagrangian subspaces $L_j(t)$ have the frames $Y_j(t)$ for $j=1,2,\dots,m$) remains constant for all $t\in [a,b]$ by \cite[Proposition A.3.8]{Kash}, then by \eqref{muthrkash} and \eqref{constdim} the cyclic sums $\mu_c^{\pm}(Y_1(t),\,Y_2(t),\dots,Y_m(t))$ and $\nu_c^{\pm}(Y_1(t),\,Y_2(t),\dots,Y_m(t))$ remain constant as well.
\end{proof}

Theorem~\ref{kashcycl} coupled with Definition~\ref{cyclsumd}, Proposition~\ref{invar}, and the index results in Section~\ref{sec2} imply the following representations of the Kashiwara index.
\begin{corollary}\label{cor1}
For the Kashiwara index \eqref{kashm} we have the following representations
\begin{equation}\label{con1-kash}
  \tau(Y_1,Y_2,\dots,Y_m)=\sum\limits_{j=1}^{m-1}\sign \mathcal P(\tilde Y_j,\tilde Y_{j+1})+\sign \mathcal P(\tilde Y_m,\tilde Y_1),\quad \tilde Y_j=R^{-1}Y_j=\binom{\tilde X_j}{\tilde U_j}
\end{equation}
where $R\in\Sp$ is arbitrary and the symmetric matrices $\mathcal P(\tilde Y_i,\tilde Y_j)$ are defined in \eqref{compind1} for the comparative indices $\mu(\tilde Y_i,\tilde Y_j),$ i.e.,
\begin{equation}\label{Pij}
  \mathcal P(\tilde Y_i,\tilde Y_j)=F_{\mathcal M}\;w^T_{i,j}\,\tilde X_i^{\dag}\,\tilde X_j\; F_{\mathcal M},\quad \mathcal M=F_{\tilde X_i}\;w_{i,j},\quad w_{i,j}:=w(Y_i,Y_j),
\end{equation}
in particular, for $R:=Z_m,$ where $Y_m=Z_m\bnn$ we have by Corollary~\ref{impinv}
\begin{equation}\label{con2-kash}
  \tau(Y_1,Y_2,\dots,Y_m)=\sum\limits_{j=1}^{m-2}\sign \mathcal P(Z_m^{-1} Y_j,Z_m^{-1}  Y_{j+1}).
\end{equation}
\end{corollary}
\begin{proof} By  Theorem~\ref{kashcycl} and Proposition~\ref{invar}
\begin{equation*}\begin{aligned}
 \tau(Y_1,Y_2,\dots,Y_m)= \mu_c^{+}(Y_1,Y_2,\dots,Y_m)&-\mu_c^{-}(Y_1,Y_2,\dots,Y_m)\\= \mu_c^{+}(R^{-1}Y_1,R^{-1}Y_2,\dots,R^{-1}Y_m)&-\mu_c^{-}(R^{-1}Y_1,R^{-1}Y_2,\dots,R^{-1}Y_m),
\end{aligned}\end{equation*}
then one can derive \eqref{con1-kash} substituting definition \eqref{cyclsum} of $\mu_c^{\pm}(R^{-1}Y_1,R^{-1}Y_2,\dots,R^{-1}Y_m)$ into the right-hand side of the last equality and using \eqref{defmu}, \eqref{defmu*} according to
\begin{equation*}\begin{aligned}
  \mu^*(\tilde Y_k,\tilde Y_{k+1})-\mu(\tilde Y_{k},\tilde Y_{k+1})&=\ind(-\mathcal P(\tilde Y_k,\tilde Y_{k+1}))-\ind(\mathcal P(\tilde Y_k,\tilde Y_{k+1}))=\sign \mathcal P(\tilde Y_k,\tilde Y_{k+1}),\,k=1,\dots,m-1,\\\mu^*(\tilde Y_m,\tilde Y_{1})-\mu(\tilde Y_{m},\tilde Y_{1})&=\sign \mathcal P(\tilde Y_m,\tilde Y_1),
\end{aligned}\end{equation*}
where $\mathcal P(\tilde Y_i,\tilde Y_j)$ are calculated according to \eqref{Pij}. The proof is completed.
\end{proof}
\begin{corollary}\label{cor2}
Under the notation of Corollary~\ref{maincon1+} and Theorem~\ref{maincon23} we have the following representations for the Kashiwara index \eqref{kashm}
\begin{equation}\label{con3-kash}
  \tau(Y_1,Y_2,\dots,Y_m)=-\sign(S_{1,2,\dots,m})=-\sign(F_{\mathcal W}\;S_{1,2,\dots,m-1}\;F_{\mathcal W})=-\sign(F_{\tilde{\mathcal M}} (S_{2,3,\dots,m-1}-\mathcal D-\mathcal D^T )F_{\tilde{\mathcal M}}).
\end{equation}
\end{corollary}
\begin{proof}
The proof follows from Theorem~\ref{kashcycl}
coupled with Theorems~\ref{maincon},~\ref{maincon23}, in particular, we apply Theorem~\ref{maincon23} cancelling the same addends $\rank \tilde{\mathcal M}$ in the representations of $\nu_c^{\pm}(Y_1,Y_2,\dots,Y_m)$ in \eqref{maincon2}.
\end{proof}
\begin{remark}
Corollary~\ref{cor1} applied to the case $m=3$ leads to the representations of the Kashiwara index
\begin{equation}\label{signcycl1}\begin{aligned}
\tau(L_1,L_2,L_3)&=\sign(\mathcal P(\tilde Y_1,\tilde Y_2))+\sign(\mathcal P(\tilde Y_2,\tilde Y_3))+\sign(\mathcal P(\tilde Y_3,\tilde Y_1))\\&=\sign(\mathcal P(Z_3^{-1} Y_1, Z_3^{-1}Y_2))=\sign(\mathcal P(Z_1^{-1} Y_2, Z_1^{-1}Y_3))=\sign(\mathcal P(Z_2^{-1} Y_3, Z_2^{-1}Y_1)),
\end{aligned}\end{equation}
where in the last row of \eqref{signcycl1} we applied Corollary~\ref{cor1} to the cases $R:=Z_3,$ $R:=Z_1$, and $R:=Z_2.$

Observe also, that the last equality in Corollary~\ref{cor2} for $m=3$ can be simplified according to Example~\ref{ex1} (see \eqref{compwron2}, \eqref{ind1})
\begin{equation}\label{signcycl2}\begin{aligned}
\tau(L_1,L_2,L_3)=\sign(F_{{\mathcal M}} \mathcal D F_{{\mathcal M}})=\sign(\mathcal P(Z_3^{-1} Y_1, Z_3^{-1}Y_2)).
\end{aligned}\end{equation}

\end{remark}
\subsection{Cyclic sums in oscillation theory of \eqref{SDS}}
In this section we consider discrete symplectic system \eqref{SDS}, where $\cS_k\in\Sp$ is separated into $n\times n$ blocks according to
\begin{equation}\label{symplm}
  \cS_k=\begin{pmatrix}
          \cA_k & \cB_k \\
          \cC_k & \cD_k \\
        \end{pmatrix}.
\end{equation}
Introduce   $2n\times n$ matrix solutions $\mathcal Y_k,\,k=0,1,\dots,N+1$ with conditions \eqref{conj} (conjoined bases of \eqref{SDS}) associated with symplectic fundamental matrices $\mathcal Z_k,\,k=0,1,\dots N+1$ such that the condition $\mathcal Y_k=\mathcal Z_k\bnn$ holds. We define the \textit{principal solution} $\mathcal Y_k^{[M]}$ of \eqref{SDS} at $k=M$ by the initial condition
\begin{equation}\label{princ}
  \mathcal Y_M^{[M]}=\bnn,\quad M=0,1,\dots,N+1.
\end{equation}

According to the definition, see \cite[Definition 1]{Kratz1} a conjoined basis $\mathcal Y_{k}=\binom{\mathcal {\mathcal X}_k}{\mathcal U_k}$ has a \textit{forward focal point} of the
multiplicity $m_1(k)$ in the point $k+1$ where $m_{1} (k) = \rank M_{k}, $
$ M_{k} = \left( {I - {\mathcal X}_{k + 1} {\mathcal X}_{k + 1}^{\dag} }  \right){\mathcal B}_{k},
$
and this basis has a forward focal point of the multiplicity $m_{2} (k) $ in
the interval $(k,k+1)$ if $m_{2} (k) = \ind(T_{k}^{T} {\mathcal X}_{k} {\mathcal X}_{k + 1}^{\dag } {\mathcal B}_{k}T_{k}) ,\;T_{k} = I -
M_{k}^{\dag}M_{k}. $ The number of forward focal points in $\left( {k,k + 1}
\right]$ is defined by $m\left( {k} \right) = m_{1} \left( {k}
\right) + m_{2} \left( {k} \right)$. This definition can be briefly rewritten in terms of the comparative index as follows (see \cite[Lemma 3.1]{Elyseeva1})
\begin{equation}\label{forfocp}
  m(k)=\mu(\mathcal Y_{k+1},\cS_k\bnn)=\mu^*(\mathcal Z_{k+1}^{-1}\bnn,\mathcal Z_{k}^{-1}\bnn).
\end{equation}
The multiplicities $m^*(k)$ of \textit{backward focal points}  associated with conjoined bases of the so-called \textit{time-reversed } symplectic system (see \cite[Section 2.1.2]{DEH})
\begin{equation}\label{timerevSDS}
  y_k=\cS^{-1}_k y_{k+1},\,k=N,N-1,\dots,0
\end{equation}
are defined as follows, see \cite[Lemma 3.2]{Elyseeva1}
\begin{equation}\label{forfocp*}
  m^*(k)=\mu^*(\mathcal Y_k,\cS_k^{-1}\bnn)=\mu(\mathcal Z_{k}^{-1}\bnn,\mathcal Z_{k+1}^{-1}\bnn).
\end{equation}
Consider the total numbers of forward and backward focal points of conjoined bases of symplectic system~\eqref{SDS}
\begin{equation}\label{forw}
  l(\mathcal Y,0,N+1)=\sum\limits_{k=0}^N \mu(\mathcal Y_{k+1},\cS_k\bnn)=\sum\limits_{k=0}^N \mu^*(\mathcal Z_{k+1}^{-1}\bnn,\mathcal Z_{k}^{-1}\bnn)
\end{equation}
and
\begin{equation}\label{backw}
  l^*(Y,0,N+1)=\sum\limits_{k=0}^N \mu^*(\mathcal Y_{k},\cS_k^{-1}\bnn)=\sum\limits_{k=0}^N \mu(\mathcal Z_{k}^{-1}\bnn,\mathcal Z_{k+1}^{-1}\bnn).
\end{equation}
For any conjoined basis $\mathcal Y_k$ of \eqref{SDS} (see \cite[formulas (4.62), (4.66)]{DEH}) the following inequalities hold
\begin{equation}\label{sepc}\begin{aligned}
 l(\mathcal Y^{[0]},0,N+1) &\le l(\mathcal Y,0,N+1)\le l(\mathcal Y^{[N+1]},0,N+1),\\
 l^*(\mathcal Y^{[N+1]},0,N+1) &\le l^*(\mathcal Y,0,N+1)\le l^*(\mathcal Y^{[0]},0,N+1),
\end{aligned}\end{equation}
moreover, by \cite[Lemma 3.3]{Elyseeva1}, \cite[Theorems 4.34, 4.35]{DEH}
\begin{equation}\label{l0=lN}
  l(\mathcal Y^{[0]},0,N+1)=l^*(\mathcal Y^{[N+1]},0,N+1),\quad  l^*(\mathcal Y^{[0]},0,N+1)=l(\mathcal Y^{[N+1]},0,N+1).
\end{equation}
Consider  cyclic sums \eqref{cyclsum}, \eqref{cyclsummod} for the special case
\begin{equation}\label{speccasc}
  m:=N+2,\,Y_k:=\mathcal Z_{k-1}^{-1}\bnn,\,k=1,\dots,N+2,
\end{equation}
where $\mathcal Z_k\in \Sp$ is a fundamental matrix of \eqref{SDS}.
The main result of this section is the following theorem.
\begin{theorem}\label{cycldisc}
The cyclic sums \eqref{secka} are invariant with respect to a choice of a fundamental matrix $\mathcal {\mathcal Z}_k\in \Sp$ of \eqref{SDS} and
\begin{equation}\label{maindiska}\begin{aligned}
\mu_c^{-}({\mathcal Z}_0^{-1}\bnn,{\mathcal Z}_1^{-1}\bnn,\dots,{\mathcal Z}_{N+1}^{-1}\bnn) &=l^*({\mathcal Y}^{[0]},0,N+1)\\=\mu_c^{+}({\mathcal Z}_{N+1}^{-1}\bnn,{\mathcal Z}_N^{-1}\bnn,\dots,{\mathcal Z}_{0}^{-1}\bnn)&=l({\mathcal Y}^{[N+1]},0,N+1).
\end{aligned}\end{equation}
By a similar way, for any choice of a symplectic fundamental matrix $\mathcal {\mathcal Z}_k$ of \eqref{SDS}
\begin{equation}\label{maindiskb}\begin{aligned}
\nu_c^{-}({\mathcal Z}_0^{-1}\bnn,{\mathcal Z}_1^{-1}\bnn,\dots,{\mathcal Z}_{N+1}^{-1}\bnn) &=l^*({\mathcal Y}^{[N+1]},0,N+1)\\=\nu_c^{+}({\mathcal Z}_{N+1}^{-1}\bnn,{\mathcal Z}_N^{-1}\bnn,\dots,{\mathcal Z}_{0}^{-1}\bnn)&=l({\mathcal Y}^{[0]},0,N+1),
\end{aligned}\end{equation}
where $l({\mathcal Y}^{[M]},0,N+1)$ and $l^*({\mathcal Y}^{[M]},0,N+1)$ are the total numbers of forward  and backward focal points of the principal solution at $M$.
\end{theorem}
\begin{proof}
 Consider \eqref{backw} for the case of the principal solution ${\mathcal Y}^{[0]}_k$ at zero, i.e., for ${\mathcal Y}^{[0]}_0=\bnn,$ and introduce the symplectic fundamental matrix ${\mathcal Z}^{[0]}_k$ such that ${\mathcal Y}^{[0]}_k={\mathcal Z}^{[0]}_k\bnn.$ We see by the second equality in \eqref{backw} and \eqref{speccase}
\begin{equation}\label{princbackw}
  l^*({\mathcal Y}^{[0]},0,N+1)=\mu_c^{-}(\bnn,{\mathcal Z}^{[0]\,-1}_1\bnn,{\mathcal Z}^{[0]\,-1}_2\bnn,\dots,{\mathcal Z}^{[0]\,-1}_{N+1}\bnn).
\end{equation}
Next, putting $R:={\mathcal Z}_0$ in Proposition~\ref{invar} we have
\begin{equation}\label{princbackwc}\begin{aligned}
\mu_c^{-}({\mathcal Z}_0^{-1}\bnn,{\mathcal Z}_1^{-1}\bnn,{\mathcal Z}_2^{-1}\bnn,&\dots,{\mathcal Z}_{N+1}^{-1}\bnn)\\ =\mu_c^{-}(\bnn,{\mathcal Z}_0{\mathcal Z}_1^{-1}\bnn,{\mathcal Z}_0{\mathcal Z}_2^{-1}\bnn,&\dots,{\mathcal Z}_0{\mathcal Z}_{N+1}^{-1}\bnn)\\
=\mu_c^{-}(\bnn,({\mathcal Z}_1{\mathcal Z}_0^{-1})^{-1}\bnn,({\mathcal Z}_2{\mathcal Z}_0^{-1})^{-1}\bnn,&\dots,({\mathcal Z}_{N+1}{\mathcal Z}_0^{-1})^{-1}\bnn)\\=\mu_c^{-}(\bnn,{\mathcal Z}^{[0]\,-1}_1\bnn,{\mathcal Z}^{[0]\,-1}_2\bnn,&\dots,{\mathcal Z}^{[0]\,-1}_{N+1}\bnn).
\end{aligned}\end{equation}
By \eqref{princbackw} and \eqref{princbackwc} we have proved the first  equality in \eqref{maindiska}. Next we use Proposition~\ref{invar1}(iii)
\begin{equation}\label{princforw}
  \mu_c^{-}({\mathcal Z}_0^{-1}\bnn,{\mathcal Z}_1^{-1}\bnn,{\mathcal Z}_2^{-1}\bnn,\dots,{\mathcal Z}_{N+1}^{-1}\bnn)=\mu_c^{+}({\mathcal Z}_{N+1}^{-1}\bnn,{\mathcal Z}_N^{-1}\bnn,\dots,{\mathcal Z}_{0}^{-1}\bnn)
\end{equation}
and consider \eqref{forw} for the principal solution ${\mathcal Y}^{[N+1]}_k$ at $N+1$, i.e., for ${\mathcal Y}^{[N+1]}_{N+1}=\bnn$  introducing the symplectic fundamental matrix ${\mathcal Z}^{[N+1]}_k$ such that ${\mathcal Y}^{[N+1]}_k={\mathcal Z}^{[N+1]}_k\bnn.$ We derive from the second equality in \eqref{forw} and \eqref{speccase}
\begin{equation}\label{princforwa}
  l({\mathcal Y}^{[N+1]},0,N+1)=\mu_c^{+}(\bnn,{\mathcal Z}^{[N+1]\,-1}_N\bnn,{\mathcal Z}^{[N+1]\,-1}_{N-1}\bnn,\dots,{\mathcal Z}^{[N+1]\,-1}_{0}\bnn)
\end{equation}
while by  invariant property \eqref{inv1}
\begin{equation}\label{invforw}\begin{aligned}
  \mu_c^+({\mathcal Z}_{N+1}^{-1}\bnn,{\mathcal Z}_N^{-1}\bnn,{\mathcal Z}_{N-1}^{-1}\bnn,&\dots,{\mathcal Z}_{0}^{-1}\bnn)\\=\mu_c^+(\bnn,{\mathcal Z}_{N+1}{\mathcal Z}_N^{-1}\bnn,{\mathcal Z}_{N+1}{\mathcal Z}_{N-1}^{-1}\bnn,&\dots,{\mathcal Z}_{N+1}{\mathcal Z}_{0}^{-1}\bnn)\\=\mu_c^+(\bnn,{\mathcal Z}_N^{[N+1]\,-1}\bnn,{\mathcal Z}_{N-1}^{[N+1]\,-1}\bnn,&\dots,{\mathcal Z}_{0}^{[N+1]\,-1}\bnn).
\end{aligned}\end{equation}
By \eqref{princforw}, \eqref{princforwa}, and \eqref{invforw} we complete the proof of \eqref{maindiska}.

For the proof of \eqref{maindiskb} we see by \eqref{backw} and \eqref{speccase}
\begin{equation}\label{princforwb}
  l^*({\mathcal Y}^{[N+1]},0,N+1)=\nu_c^{-}({\mathcal Z}^{[N+1]\,-1}_0\bnn,{\mathcal Z}^{[N+1]\,-1}_{1}\bnn,\dots,{\mathcal Z}^{[N+1]\,-1}_{N}\bnn,\bnn),
\end{equation}
and by putting $R:={\mathcal Z}_{N+1}$ in Proposition~\ref{invar}
\begin{equation}\label{diskb}\begin{aligned}
\nu_c^{-}({\mathcal Z}_0^{-1}\bnn,{\mathcal Z}_1^{-1}\bnn,&\dots,{\mathcal Z}_{N}^{-1}\bnn,{\mathcal Z}_{N+1}^{-1}\bnn) \\=\nu_c^{-}({\mathcal Z}_{N+1}{\mathcal Z}_0^{-1}\bnn,{\mathcal Z}_{N+1}{\mathcal Z}_1^{-1}\bnn,&\dots,{\mathcal Z}_{N+1}{\mathcal Z}_{N}^{-1}\bnn,\bnn)\\
=\nu_c^{-}(({\mathcal Z}_0{\mathcal Z}_{N+1}^{-1})^{-1}\bnn,({\mathcal Z}_1{\mathcal Z}_{N+1}^{-1})^{-1}\bnn,&\dots,({\mathcal Z}_NZ_{N+1}^{-1})^{-1}\bnn,\bnn)\\=\nu_c^{-}({\mathcal Z}^{[N+1]\,-1}_0\bnn,{\mathcal Z}^{[N+1]\,-1}_{1}\bnn,&\dots,{\mathcal Z}^{[N+1]\,-1}_{N}\bnn,\bnn) \end{aligned}\end{equation}
By \eqref{princforwb} and \eqref{diskb} we  prove the first equality in \eqref{maindiskb} and  then again use Proposition~\ref{invar1}(iii)
\begin{equation}\label{diska}
  \nu_c^{-}({\mathcal Z}_0^{-1}\bnn,{\mathcal Z}_1^{-1}\bnn,\dots,{\mathcal Z}_{N+1}^{-1}\bnn)=\nu_c^{+}({\mathcal Z}_{N+1}^{-1}\bnn,{\mathcal Z}_N^{-1}\bnn,\dots,{\mathcal Z}_{0}^{-1}\bnn).
\end{equation}
Finally we apply the invariant property to the right-hand side of \eqref{diska} for $R:={\mathcal Z}_0$ and prove
\begin{equation}\label{diskc}\begin{aligned}
  \nu_c^{+}({\mathcal Z}_{N+1}^{-1}\bnn,{\mathcal Z}_N^{-1}\bnn,\dots,{\mathcal Z}_{0}^{-1}\bnn)=\nu_c^{+}({\mathcal Z}_0{\mathcal Z}_{N+1}^{-1}\bnn,{\mathcal Z}_0{\mathcal Z}_N^{-1}\bnn,&\dots,{\mathcal Z}_0{\mathcal Z}_{1}^{-1}\bnn,\bnn)\\=
  \nu_c^{+}({\mathcal Z}_{N+1}^{[0]\,-1}\bnn,{\mathcal Z}_N^{[0]\,-1}\bnn,\dots,{\mathcal Z}_{1}^{[0]\,-1}\bnn,\bnn)&=l({\mathcal Y}^{[0]},0,N+1)
\end{aligned}\end{equation}
according to \eqref{forw} for ${\mathcal Y}_k:={\mathcal Y}^{[0]}_k.$ By \eqref{diska}, \eqref{diskc} the proof of \eqref{maindiskb} is completed.
\end{proof}
\begin{remark}
In the proof of Theorem~\ref{cycldisc} we used only  definitions \eqref{forw}, \eqref{backw} of the numbers of forward and backward focal points as well as Propositions~\ref{invar},~\ref{invar1}(iii) presenting the new proof of identities \eqref{l0=lN}. Moreover, applying \eqref{cyclsum}, \eqref{cyclsummod} to the cyclic sums in Theorem~\ref{cycldisc} one can also prove inequalities \eqref{sepc} based on separation results in \cite[Corollary 3.1, formulas (3.9),(3.10)]{Elyseeva1}, \cite[Section 4.2.3]{DEH}. For example, by \eqref{maindiskb},\eqref{cyclsummod}, and \eqref{forw}
\begin{equation*}
  l({\mathcal Y}^{[0]},0,N+1)=\nu_c^{+}({\mathcal Z}_{N+1}^{-1}\bnn,{\mathcal Z}_N^{-1}\bnn,\dots,{\mathcal Z}_{0}^{-1}\bnn)=l({\mathcal Y},0,N+1)-\mu^*({\mathcal Z}_{N+1}^{-1}\bnn,{\mathcal Z}_{0}^{-1}\bnn)\ge 0,
\end{equation*}
where by Lemma~\ref{prop ind}(iv) $\mu^*({\mathcal Z}_{N+1}^{-1}\bnn,{\mathcal Z}_{0}^{-1}\bnn)=\mu({\mathcal Y}_{N+1},{\mathcal Z}_{N+1}{\mathcal Z}_{0}^{-1}\bnn)=\mu({\mathcal Y}_{N+1},{\mathcal Y}_{N+1}^{[0]}),$ compare with \cite[Corollary 3.1]{Elyseeva1}. Observe that we  also proved the lower bound in \eqref{sepc}. By a similar way one can  derive other equalities in \cite{Elyseeva1}, \cite[Section 4.2.3]{DEH}.
\end{remark}

Applying Theorems~\ref{maincon},~\ref{maincon23} we derive the following representation for the numbers of focal points  in terms of the indices of symmetric  matrices.
\begin{theorem}\label{indfocd}
Let $\mathcal Y_k^{[M]}$ be the \textit{principal solutions}  of \eqref{SDS} at $k=M$ for $M=0,1,\dots,N+1$ with the upper blocks $\mathcal X_k^{[M]}.$ Then  we have the following representations for the number of focal points of \eqref{SDS} given by \eqref{maindiska}, \eqref{maindiskb}
\begin{equation}\label{princbackwind}\begin{aligned}
  l^*(\mathcal Y^{[0]},0,N+1)&=l({\mathcal Y}^{[N+1]},0,N+1)=\ind(-S_{1,2,\dots,N+2}^{[0]}),\;N\ge 0,\\
   l(\mathcal Y^{[0]},0,N+1)&=l^*({\mathcal Y}^{[N+1]},0,N+1)=\ind(-\bar S_{1,2,\dots,N+1}^{[0]}),\;N\ge 1,\\
   \bar S_{1,2,\dots,N+1}^{[0]}&=\begin{pmatrix}
                                   0 & \tilde{\mathcal M_d} \\
                                   \tilde{\mathcal M_d}^T & S_{2,3,\dots,N+1}^{[0]}-\mathcal D_d-\mathcal D_d^T\\
                                 \end{pmatrix},
\end{aligned}\end{equation}
where $S_{2,3,\dots,N+1}^{[0]}$ is the submatrix of $S_{1,2,\dots,N+2}^{[0]}.$
The $n\times n$ blocks  $S_{1,2,\dots,N+2}^{[0]}(i,j),$  of $S_{1,2,\dots,N+2}^{[0]}$ are defined as follows
  \begin{equation}\label{princbackwindel}\begin{aligned}
  S_{1,2,\dots,N+2}^{[0]}(i,j)&=\mathcal X^{[i-1]\,T}_{j-1}=-\mathcal X^{[j-1]}_{i-1},\quad S_{1,2,\dots,N+2}^{[0]}(j,i)=S_{1,2,\dots,N+2}^{[0]\,T}(i,j),\;j\ge i,\quad i,j=1,2,\dots,N+2,\\
 \end{aligned}\end{equation}
in particular,
\begin{equation}\label{seconddiag}
\begin{aligned}
  S_{1,2,\dots,N+2}^{[0]}(i+1,i)&=S_{1,2,\dots,N+2}^{[0]\;T}(i,i+1)=\mathcal B_{i-1},
 \end{aligned}\end{equation}
where $\mathcal B_{k}$ is the block of $\cS_k$ according to \eqref{symplm}.
We also have
\begin{equation}\label{MKNd}\begin{aligned}
\tilde{\mathcal M_d}&=(I-\mathcal X_0^{[N+1]\;\dag}\mathcal X_0^{[N+1]})\mathcal N_d,\quad
\mathcal D_d=-\mathcal K^T_d\mathcal X_{N+1}^{[0]\,\dag}\mathcal N_d,\\
\mathcal K_d&=(\mathcal X_1^{[0]\;T}\;\mathcal X_2^{[0]\;T} \dots \mathcal X_N^{[0]\;T}),\quad
\mathcal N_d=(\mathcal X_1^{[N+1]\;T}\;\mathcal X_2^{[N+1]\;T} \dots \mathcal X_N^{[N+1]\;T}),
   \end{aligned}\end{equation}
where $\tilde{\mathcal M_d}$ in \eqref{princbackwind} can be replaces by
\begin{equation}\label{Md}
  {\mathcal M_d}=(I-\mathcal X_{N+1}^{[0]\;\dag}\mathcal X_{N+1}^{[0]})\mathcal K_d.
\end{equation}
\end{theorem}
\begin{proof}
The first equalities in  \eqref{princbackwind} follow from \eqref{maindiska} coupled with Theorem~\ref{maincon}, where we derive \eqref{princbackwindel} according to \eqref{Sm}. Indeed,
by  \eqref{maindiska}  we have $$l^*(\mathcal Y^{[0]},0,N+1)=l(\mathcal Y^{[N+1]},0,N+1)=\mu_c^{-}(\mathcal Z^{-1}_0\bnn,\mathcal Z^{-1}_1\bnn,\dots,\mathcal Z^{-1}_{N+1}\bnn) ,$$ then by \eqref{Sm} applied to the case \eqref{speccasc}  we derive
 \begin{equation}\label{add11}\begin{aligned}
  S_{1,2,\dots m}^{[0]}(i,j)&=w(\mathcal Z^{-1}_{i-1}\bnn,\mathcal Z^{-1}_{j-1}\bnn)=w(\mathcal Z_{j-1}\mathcal Z^{-1}_{i-1}\bnn,\bnn)=\mathcal X_{j-1}^{[i-1]\;T}\\&=w(\bnn,\mathcal Z_{i-1}\mathcal Z^{-1}_{j-1}\bnn)=-\mathcal X_{i-1}^{[j-1]},\,j\ge i.
\end{aligned}\end{equation}

Note that $S_{1,2,\dots,N+2}^{[0]}(i,i)=\mathcal X_{i-1}^{[i-1]}=0$
according to  definition \eqref{princ} of the principal solution at the point $M=i-1$
 and by \eqref{princ}, \eqref{SDS}  $\mathcal X_{i}^{[i-1]}=\mathcal B_{i-1},\,i=1,2,\dots N+1.$

The  equalities in the second row of  \eqref{princbackwind} are derived using \eqref{maindiskb}, Corollary~\ref{maincon1+}, and Theorem~\ref{maincon23}. In more details, we have $l(\mathcal Y^{[0]},0,N+1)=l^*({\mathcal Y}^{[N+1]},0,N+1)=\nu_c^{-}({\mathcal Z}_0^{-1}\bnn,{\mathcal Z}_1^{-1}\bnn,\dots,{\mathcal Z}_{N+1}^{-1}\bnn)$ by \eqref{maindiskb} and then one can calculate $\nu_c^{-}({\mathcal Z}_0^{-1}\bnn,{\mathcal Z}_1^{-1}\bnn,\dots,{\mathcal Z}_{N+1}^{-1}\bnn)$ according to the second equality in \eqref{maincon2} using the definition of $\tilde{\mathcal M_d},\,{\mathcal M_d},\,\mathcal K_d,\,\mathcal N_d$ and $\mathcal D_d$ in Corollary~\ref{maincon1+}. Observe also that in these computations we again used connection \eqref{add11}    for the upper blocks of the principal solutions $\mathcal Y_k^{[i-1]},\,\mathcal Y_k^{[j-1]}$ of \eqref{SDS}. The proof is completed.
\end{proof}
Recall that according to \cite[Theorems 1,2]{BDosly},\cite[Theorems 2.36, 2.41]{DEH} systems \eqref{SDS} and \eqref{timerevSDS}   are {disconjugate} on $[0,N+1]$ if and only if $l(\mathcal Y^{[0]},0,N+1)=l^*(\mathcal Y^{[N+1]},0,N+1)=0.$  Based on Theorem~\ref{indfocd} one can formulate the following criterion for disconjugacy of \eqref{SDS} and \eqref{timerevSDS}.
\begin{corollary}
Systems \eqref{SDS} and \eqref{timerevSDS} are disconjugate on $[0,N+1]$ if and only if the matrix  $\bar S_{1,2,\dots,N+1}^{[0]}$ defined in \eqref{princbackwind} is nonpositive definite, i.e. $\bar S_{1,2,\dots,N+1}^{[0]}\le 0.$ The last condition is equivalent to
\begin{equation}\label{nonpositive}
  \tilde{\mathcal M_d}={\mathcal M_d}=0,\quad S_{2,\dots,N+1}^{[0]}-\mathcal D_d-\mathcal D_d^T\le 0,
\end{equation}
where we use the notation in Theorem~\ref{indfocd}.
\end{corollary}
\begin{proof}
We have by Theorem~\ref{indfocd} that the condition $l(\mathcal Y^{[0]},0,N+1)=l^*(\mathcal Y^{[N+1]},0,N+1)=0$ is equivalent to $\ind(-\bar S_{1,2,\dots,N+1}^{[0]})=0.$ Then the first claim of this corollary is proved, while conditions \eqref{nonpositive} follow from  the last equality in \eqref{maincon2}.
\end{proof}
\section*{Acknowledgments}
This research is supported by  the Ministry of Science and Higher Education of the Russian Federation under project
0707-2020-0034 and by the Czech Science Foundation under grant
GA19--01246S.

\end{document}